\documentclass[a4paper, 11pt]{article}

\usepackage[margin=1in]{geometry}
\usepackage{fullpage}

\usepackage{amsmath} \allowdisplaybreaks
\usepackage{amssymb}
\usepackage{mathtools}

\usepackage{natbib}
\usepackage{bibentry}

\usepackage{algorithm}
\usepackage{algpseudocode}
\usepackage{multirow}
\usepackage{booktabs} 
\usepackage{mathrsfs}
\usepackage{xcolor}
\usepackage{bbm}
\usepackage{tikz}
\usepackage{subcaption}
\usetikzlibrary{positioning, calc}
\usepackage{etoolbox} 
\makeatletter
\def\do#1{\@namedef{#1c}{\ensuremath{\mathcal{#1}}}}
\docsvlist{A,B,C,D,E,F,G,H,I,J,K,L,M,N,O,P,Q,R,S,T,U,V,W,X,Y,Z}
\makeatother
\usepackage{colortbl}
\newcommand{\RNum}[1]{\uppercase\expandafter{\romannumeral #1\relax}}
\usepackage{makecell}
\usepackage{hyperref}

\usepackage{authblk}

\author[a]{Sasan Mahmoudinazlou}
\author[a]{Abhay Sobhanan}
\author[a]{Hadi Charkhgard}
\author[b]{Ali Eshragh}
\author[c]{George Dunn}
\affil[a]{Department of Industrial and Management Systems Engineering, University of South Florida, Tampa, FL, USA}
\affil[b]{Carey Business School, Johns Hopkins University, Washington, D.C., USA
}
\affil[c]{School of Information and Physical Sciences, University of Newcastle, NSW, Australia
}

\title{Deep Reinforcement Learning for Dynamic Order Picking in Warehouse Operations} 

\date{}

\begin{document}

\maketitle

\begin{abstract}
Order picking is a pivotal operation in warehouses that directly impacts overall efficiency and profitability.
This study addresses the dynamic order picking problem, a significant concern in modern warehouse management, where real-time adaptation to fluctuating order arrivals and efficient picker routing are crucial. 
Traditional methods, which often depend on static optimization algorithms designed around fixed order sets for the picker routing, fall short in addressing the challenges of this dynamic environment.
To overcome these challenges, we propose a Deep Reinforcement Learning (DRL) framework tailored for single-block warehouses equipped with an autonomous picking device. 
By dynamically optimizing picker routes, our approach significantly reduces order throughput times and unfulfilled orders, particularly under high order arrival rates. 
We benchmark our DRL model against established algorithms, utilizing instances generated based on standard practices in the order picking literature. Experimental results demonstrate the superiority of our DRL model over benchmark algorithms. For example, at a high order arrival rate of 0.09 (i.e., 9 orders per 100 units of time on average), our approach achieves
an order fulfillment rate of approximately 98\%, compared to the 82\% fulfillment rate observed with benchmarking algorithms.
We further investigate the integration of a hyperparameter in the reward function that allows for flexible balancing between distance traveled and order completion time. 
Finally, we demonstrate the robustness of our DRL model on out-of-sample test instances.
\end{abstract}   


\section{Introduction} \label{sec:intro}

Order picking, the process of retrieving specific items from warehouse storage to fulfill customer orders, is a labor-intensive and time-consuming operation. 
Order pickers traverse the warehouse, collecting items from various locations before delivering them to a designated depot.
Notably, order picking can constitute about 50\% of the total order processing time \citep{charkhgard2015efficient, silva2020integrating}. 
Recent growth in e-commerce, coupled with heightened customer expectations for fast delivery and increased market competition, has escalated the demand for rapid order fulfillment \citep{marchet2015investigating, rasmi2022wave,van2018designing}.
Efficient order picking, which necessitates solving a routing problem, is crucial for several reasons.
It directly enhances operational efficiency, leading to increased order throughput within the warehouse and improved customer satisfaction \citep{giannikas2017interventionist}. 
Additionally, streamlining the picking process translates to cost reductions for the business.

In modern warehouse operations, the ability to accurately forecast demand and dynamically adjust order picking decisions is paramount \citep{dauod2022real,lu2016algorithm}.
However, the existing literature predominantly focuses on exact and heuristic methods that assume a fixed set of orders is known in advance.
While this approach is suitable for static scenarios, it often fails to effectively address the dynamic nature of real-world warehouse operations. 
To overcome these challenges, a dynamic and adaptive approach to demand forecasting and picker routing is essential. 
This would enable warehouses to respond swiftly to real-time order arrivals, optimizing routing decisions, minimizing operational costs, and maximizing overall efficiency.

This research investigates picker-to-parts systems within single-block warehouses, focusing on how an autonomous picker navigates the aisles to retrieve items. This focus aligns with the growing trend of warehouse automation \citep{bansal2021stochastic, mirzaei2021impact, roy2019robot, winkelhaus2021towards, winkelhaus2022hybrid} 
and assumes that the picking device flawlessly executes the recommended routing decisions, thereby eliminating behavioral factors associated with human pickers.
To address the order picking problem, we employ a Deep Reinforcement Learning (DRL) framework, well-suited for its ability to handle sequential decision-making under uncertainty.
This approach is particularly well-suited for the real-world warehouse environment \citep{begnardi2024deep,CHARKHGARD2025DRLBOKP}, characterized by the constantly changing nature of order arrivals, item locations, picker location, and picker availability.
By leveraging DRL's ability to learn and adapt to the environment while dynamically forecasting order arrivals and optimizing picker routes, we can achieve significant improvements in warehouse efficiency and throughput.
Additionally, we show that traditional deep neural network architectures are sufficient to meet these goals, avoiding the complexity, longer training times, and higher computational demands of state-of-the-art DRL models. 

Our main contribution is a DRL framework that integrates order assignment and picker routing within a unified learning environment. By incorporating implicit batching strategies and demand forecasting, this approach effectively optimizes dynamic picker routing in single-block warehouses.
We demonstrate that our dynamic order-picking method significantly reduces order completion times and provide benchmarking against established algorithms from the literature.
Our study serves as a foundation for future research to explore the complexities of real-time order picking in multi-block warehouses with multiple coordinated pickers.

Our experiments demonstrate that as order arrival rates increase, relying solely on the shortest-distance for order picking becomes inefficient.
In scenarios with high order frequencies, the policies learned by our DRL agent substantially enhances the order fulfillment efficiency.
For example, when orders follow a Poisson process with an arrival rate ($\lambda$) of 0.09, our approach reduces order throughput time by a substantial 420\% compared to existing benchmark algorithms. 
Our approach maintains an unfulfilled order rate of approximately 2\% within the work shift, while the benchmark algorithms result in up to 18\% of orders remaining unfulfilled.
We also incorporate a hyperparameter ($\alpha$) into the reward function, enabling effective balancing of the trade-off between distance traveled and order throughput times, thus aligning with the preferences of decision-makers.

The remainder of this paper is organized as follows. 
Section~\ref{sec:lit_rev} reviews relevant literature, examining various optimization approaches for order picking and applications of deep reinforcement learning in warehouse operations. 
In Section~\ref{sec:problem}, we formally describe the problem addressed in this study. 
Section~\ref{sec:method} details our proposed solution methodology and the neural network architecture employed.
Section~\ref{sec:comp_exp} presents a series of computational experiments designed to evaluate the effectiveness of our approach against existing benchmark algorithms. 
We also analyze the robustness of our model when confronted with out-of-sample instances. 
Finally, Section~\ref{sec:conclusions} summarizes our findings with concluding remarks, and proposes potential avenues for future research in this domain.

\section{Literature Review} \label{sec:lit_rev}

In this section, we review the literature on traditional optimization methods and DRL applications for order batching and picker routing in warehouses, emphasizing their influence on our research design and methodological choices. 
Traditional optimization techniques, including exact and heuristic methods, form the foundation upon which our problem definition and benchmarks are based. Recent DRL-based approaches in warehouse operations highlight the potential for learning-based solutions in dynamic environments, directly motivating our adoption of DRL as a solution method. 
Additionally, we briefly explore DRL applications in related logistics domains characterized by dynamic task allocation and routing, reinforcing the suitability of DRL for addressing dynamic order picking challenges. 
Finally, we identify a significant research gap regarding real-time integration of dynamic batching and picker routing, situating our proposed DRL-based framework as a meaningful advancement in the field.

\subsection{Optimization Approaches for Order Picking}
Numerous optimization techniques have been developed in the literature for order picking, but these generally address static problem settings. A notable contribution is by \citet{ratliff1983order}, who designed an efficient algorithm that leverages the warehouse graph structure to optimize single-picker routing in single-block layouts. Later adaptations of this approach accommodate diverse warehouse configurations, such as two-block \citep{roodbergen2001routing}, fish-bone \citep{ccelk2014order}, and chevron layouts \citep{masae2020optimal}. Alternatively, \citet{scholz2016new} and \citet{pansart2018exact} offer mathematical models that outperform standard Traveling Salesperson Problem (TSP) solvers by eliminating redundant routes, leveraging the warehouse structure to optimize routing. However, these static methods lack adaptability to dynamic order arrivals, a key characteristic of real-world applications.

Addressing the challenges of dynamic warehouse environments, \citet{lu2016algorithm} introduced a first-come-first-served order batching policy that adapts the static algorithm to handle dynamic order inflows. Their experiments indicate notable reductions in average order completion times compared to both static exact algorithm and traditional heuristics. This adaptation is a precedent for dynamic routing solutions, and forms a baseline for our proposed DRL method to further improve efficiency in dynamic warehouse settings.

Heuristic algorithms have received significant attention due to their simplicity and ease of implementation. The S-shape heuristic, for example, offers a straightforward routing policy by assigning picker routes through aisles in an S-shaped curve. \citet{hall1993distance} presents a comparative analysis of heuristic methods, showing the effectiveness of the largest-gap heuristic, particularly with smaller pick lists. Meta-heuristic algorithms have also been explored in the literature, though they often require careful adaptation and specific local search strategies to suit different warehouse designs. 
For instance, \citet{schrotenboer2017order} presents a hybrid genetic algorithm to determine routes that simultaneously handle the pickup of products in response to customer orders and delivery of returned products to storage locations. 
However, the nature of heuristic algorithms can lead to suboptimal outcomes, which may reduce overall system efficiency \citep{roodbergen2001layout}. Notably, \citet{lu2016algorithm} demonstrated that their dynamic programming algorithm outperforms the largest-gap heuristic in terms of both order completion time and total travel distance, justifying our choice to benchmark our model against \citet{lu2016algorithm}.

Order batching is another critical strategy that can significantly reduce the distance traveled by order pickers \citep{de2007design, henn2012algorithms, muter2022order}.
\citet{menendez2017variable} employ a Variable Neighborhood Search (VNS) to group warehouse orders into picker-managed batches.  
However, most studies treat batching and routing as separate problems, simplifying one to optimize the other. 
An exception is \citet{scholz2017order}, who propose an integrated approach that combines batching with picker routing, demonstrating superior results. 
In \citet{aerts2021joint}, order batching and routing are tackled as a static vehicle routing problem, and designs a VNS heuristic for the warehouse environment. 
Recent research has seen a growing interest in dynamic pick list updating, notably in \citet{dauod2022real} and \citet{yang2021flow}.
Despite these advances, heuristics for dynamic scenarios still pose challenges, as the iterative nature of neighborhood searches can lead to computational delays. Our method leverages DRL for real-time decision-making, effectively eliminating cumulative wait times often caused by the computational time associated with optimization approaches.
For a more comprehensive literature review on order picker routing and warehouse operational procedures, readers may refer to \citet{de2007design}, \citet{van2018designing}, and \citet{masae2020order}.

\subsection{Deep Reinforcement Learning for Warehouse Operations}
Deep Reinforcement Learning (DRL) has garnered increasing attention in the warehousing domain, with recent studies primarily focusing on order batching and assignment problems. 
\citet{cals2021solving} pioneered this approach, using DRL to minimize tardy orders in order batching. Subsequently, \citet{beeks2022deep} extended this method to a multi-objective problem, considering picking costs and analyzing the trade-off between order tardiness and picker efficiency. 
Similarly, \citet{li2019task} employed DRL for assigning tasks to autonomous robots, aiming to minimize makespan in the event of traffic conflicts. 
Some recent studies also utilize DRL for warehouse inventory management \citep{kaynov2024deep, tian2024iacppo}. 
These studies demonstrate the potential of DRL, yet they primarily address order batching rather than the intricacies of dynamic picker routing by minimizing travel distances or order completion times.

In picker routing, \citet{dunn2024deep} introduced a DRL model for the order picking as a combinatorial optimization problem. The authors employ an attention-based neural network that generates high-quality solutions and proves particularly effective for large-scale problems where exact methods become computationally expensive. However, this approach is limited to static scenarios. 
Echoing the need for real-world applications, \citet{roodbergen2009survey} advocates dynamic models, underscoring the relevance of our research work. More recently, \citet{begnardi2024deep} applied graph neural networks in DRL for collaborative human-robot picking, formulated as an online bipartite matching problem. However, their evaluation is limited to grid-based warehouse layouts with several pickers, and benchmarks are solely variations of the greedy heuristic. 
Given the large number of pickers within their small-scale test instances, the study places greater emphasis on the order assignment problem rather than the routing problem.

We would like to acknowledge other recent studies in the literature, such as \citet{cheng2024deep} and \citet{neves2024learning}, which illustrate the flexibility of DRL in targeted warehouse applications, including human–picker interactions and autonomous shelving transport. While these works underscore the potential of DRL in warehouse environments, their focus and complexity differ significantly from the dynamic order batching and picker routing problem addressed in our study. Collectively, these pioneering efforts (including those discussed earlier in this section) affirm the suitability of DRL for dynamic and complex decision-making in warehouse environments, thereby supporting our methodological choice and underscoring the distinct contribution of our work.

\subsection{Broader DRL Applications for Dynamic Assignment and Routing}

Our proposed model for order picking aligns with broader applications of DRL in dynamic task allocation and routing, where tasks arrive unpredictably and must be assigned in real time. We highlight here several related studies that, while addressing distinct problem domains, share methodological similarities with our work. 
For example, \citet{da2023policies} investigate the traveling maintainer problem, in which a single maintainer responds to dynamically emerging maintenance requests. As with our order picker handling incoming orders, this involves real-time routing, task prioritization, and decision-making under uncertainty. Both studies use DRL to optimize long-term efficiency; however, their objective is to minimize maintenance costs for dispersed assets, whereas we focus on reducing order wait times within a structured warehouse layout.

DRL has also been applied to various other transportation and logistics challenges. \citet{ying2020actor} use an actor-critic DRL approach for metro train scheduling under stochastic demand, demonstrating DRL’s effectiveness in large-scale dynamic control. 
\citet{akkerman2025comparison} compare RL-based value function approximation and policy function approximation for dynamic vehicle routing with stochastic customer requests, adjusting routes in response to new demand. Similarly, \citet{kitchat2024deep} develop a Deep Deterministic Policy Gradient algorithm for distributing epidemic prevention materials during public health emergencies, illustrating the adaptability of DRL to real-time allocation and routing decisions in supply chain logistics. 
These studies underscore DRL’s versatility in transportation logistics, particularly in vehicle routing and scheduling. In contrast, our approach applies DRL to warehouse order picking, where efficiency gains emerge from optimizing picker movement and batching strategies.

In inventory management, \citet{dehaybe2024deep} apply DRL to lot-sizing under non-stationary demand, demonstrating the method’s ability to learn near-optimal replenishment policies without frequent re-optimization. Similarly, \citet{oroojlooyjadid2022deep} employ a Deep Q-Network in the Beer Game, showing how DRL can optimize decentralized, multi-stage ordering decisions even under irrational agent behavior. By outperforming base-stock policies and adapting rapidly to changing conditions, their approach illustrates DRL’s flexibility in complex inventory scenarios. Together, these examples underscore DRL’s versatility across diverse supply chain and logistics challenges.

Despite these advancements, existing studies do not provide a comprehensive DRL framework tailored specifically to the warehouse order picking problem. Our contribution lies in modeling warehouse-specific constraints, real-time batching, and picker routing strategies. Further, our reward function explicitly captures key efficiency metrics in warehouse operations, ensuring practical applicability.

\color{black}
\subsection{Research Gap}

\begin{table}[ht]
    \caption{Comparison of Related Order Picking Literature}
    \label{tab:literature}
    \resizebox{\textwidth}{!}{
    \begin{tabular}{@{}lcccccm{7cm}@{}}
        \toprule
        \multicolumn{1}{c}{\multirow{3}{*}{Reference}} & 
        \multicolumn{1}{c}{\multirow{2}{*}{Order}} & 
        \multicolumn{1}{c}{\multirow{2}{*}{Online}} & 
        \multicolumn{1}{c}{\multirow{1}{*}{Single-Block}} & 
        \multicolumn{1}{c}{\multirow{2}{*}{Minimization}} & 
        \multicolumn{1}{c}{\multirow{3}{*}{Methodology}} & 
        \multicolumn{1}{c}{\multirow{3}{*}{Comments}} \\
        & \multicolumn{1}{c}{\multirow{2}{*}{Batching}} & 
        \multicolumn{1}{c}{\multirow{2}{*}{Routing}} & 
        \multicolumn{1}{c}{\multirow{1}{*}{Warehouse}} & 
        \multicolumn{1}{c}{\multirow{2}{*}{Objective}} & & \\  
        & & & \multicolumn{1}{c}{Layout} & & & \\
        \midrule        
    \citet{ratliff1983order}   &        &               & \checkmark                 &  travel distance               & Exact                           & Solves the static picker routing problem.                        \\
    \citet{lu2016algorithm}              &        & \checkmark          & \checkmark                 &  travel distance               & Exact                           & Orders are batched on a first-come, first-served basis.             \\
    \citet{scholz2016new}                    & \checkmark   &               & \checkmark                 &  travel distance               & \makecell{Exact+ \\ Heuristic}                 & Uses Iterated Local Search for order batching, combined with an exact algorithm for routing optimization.          \\
    \citet{aerts2021joint}                    & \checkmark   &               & \checkmark                 &  travel distance               & Heuristic                       & Models the problem as a VRP variant and applies Variable Neighborhood Search to find solutions.          \\
    \citet{cals2021solving}            & \checkmark   &               & \checkmark                 &  tardy orders                  & DRL & Order sequencing for routing is determined using a heuristic approach.                     \\
    \citet{beeks2022deep}                    & \checkmark   &               & \checkmark                 &  \makecell{tardy orders \& \\ picking cost} & DRL & Extends \citet{cals2021solving} to incorporate two objectives.                       \\
    \citet{dunn2024deep}            &        &               & \checkmark                 & travel distance               & DRL & Solves the static picker routing problem using DRL.           \\
    \citet{begnardi2024deep}        & \checkmark   & \checkmark          &                      & order wait time               & DRL & Grid-based warehouse layout; 
    Focuses on order assignment to multiple pickers to reduce delays; 
    Baselines lack TSP routing or traditional order picking heuristics. \\
    \hline
    This Study                    & \checkmark   & \checkmark          & \checkmark                 & \makecell{order wait time \& \\ travel distance}               & DRL &    --      \\ \bottomrule          
    \end{tabular}
    }
    \end{table}

Our work is the first comprehensive investigation into dynamic order picker routing using DRL, allowing the agent to implicitly learn order batching strategies within a unified learning environment. 
Table~\ref{tab:literature} provides a comparative analysis of related order picking literature, highlighting the unique contributions of our study.
By integrating assignment and routing, our DRL agent provides rapid, on-the-fly decision-making, crucial for reducing delays in dynamic settings. We prioritize the development of a simple and effective model for single-block warehouses, showing its efficacy in minimizing order processing times, thereby addressing a clear gap in existing research.

\section{Problem Description} \label{sec:problem}

\begin{figure}
    \centering
    \includegraphics[width=0.8\textwidth]{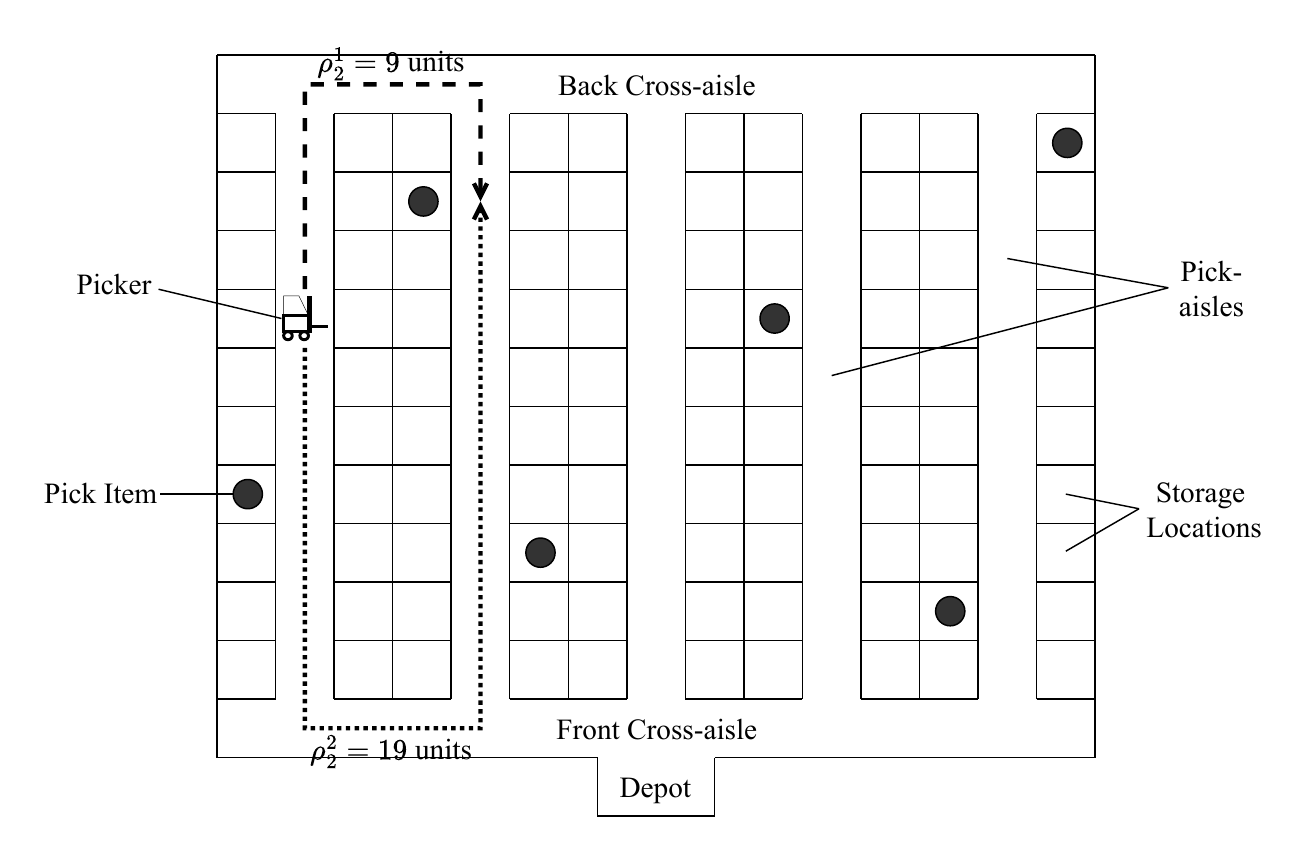}
    \caption{An illustrative example of warehouse layout for the order-picking problem}
    \label{fig:warehouse}
\end{figure}

Consider the order-picking problem for a single picker within a rectangular warehouse. The warehouse consists of a single block with multiple parallel aisles, interconnected by two cross-aisles located at the ends. Inventory storage locations are positioned on both sides of each aisle.
Figure~\ref{fig:warehouse} provides a visual representation of this warehouse layout. The picker must retrieve requested items from designated storage locations within the warehouse. Movement is restricted to the aisles, where the picker can access storage locations on either side. However, there is no lateral movement within an aisle.
To transition between aisles, the picker must use one of the two cross-aisles. These cross-aisles do not contain any storage locations and serve solely as pathways for aisle transitions. Notably, the picker can only enter or exit an aisle through one of these cross-aisles.

The order-picking device or vehicle utilized for transporting these items has a total capacity of $K$ items.
Initially, the picker starts from the depot and moves through the warehouse to collect the assigned items and delivers them to the depot, subject to the capacity constraint at all times. Figure~\ref{fig:warehouse} illustrates that this warehouse layout offers two alternative paths the picker can take to collect the next assigned order, denoted by the distances $\rho$.  We assume that multiple orders may be received from across the warehouse at any time, with each order representing a single unit and having the same level of priority. Additionally, we assume the requested items are identical in size, a common assumption in the literature \citep{aerts2021joint,lu2016algorithm,pansart2018exact, scholz2016new}. This assumption is not overly restrictive, as the smallest item in the warehouse can serve as a reference unit (e.g., 1 unit), with larger items represented as multiples of this unit to account for size differences. To prevent fractional item pickups in such scenarios, the proposed method (presented in the next section) can handle this seamlessly without requiring re-training. Specifically, orders can be temporarily excluded from the system if the picker's capacity is insufficient to pick up the entire order. These excluded orders are reintegrated into the system once the picker completes a drop-off at the depot.

We will discuss the notation details of distances in the next section. The primary objective is to minimize the average waiting time for orders, and thereby increasing the throughput in warehouse operations. 
The picker route is determined based on the demand forecast, and a dynamically updated pick-list, which is adjusted based on incoming orders and the remaining orders to be collected. 
This optimization problem is considered over a time horizon that resembles a human picker's work shift.

A key distinction of our problem setting, compared to the majority of the literature, is that it addresses a dynamic order-picking problem rather than a static one, where predefined pick lists are optimized. In our setting, orders arrive dynamically into the system, and demand fluctuates throughout the work shift.  This dynamic environment necessitates real-time decision-making for order batching and traversal routes to minimize order wait times. Reducing these wait times not only enhances customer satisfaction but also frees up resources to fulfill additional orders—an especially critical factor in the rapidly expanding e-commerce sector.  Our approach further enhances order picking by incorporating \textit{implicit demand forecasting}, enabling pickers to anticipate future orders and optimize their paths proactively. This means that instead of taking explicit forecasting steps, our DRL agent implicitly integrates forecasting into its actions. Specifically, the DRL agent's decisions reflect an understanding of potential future orders, even if they do not align directly with current unfulfilled orders.  For example, the agent may opt to wait instead of taking immediate action, despite pending orders in the system. Alternatively, it may navigate to aisles or locations that are currently empty, anticipating that future orders will require those areas. This behavior contrasts sharply with traditional methods, which rely solely on observed order arrivals and re-optimize the system each time a new order appears, without accounting for future demand.   The effectiveness of this indirect approach to demand forecasting is demonstrated in our numerical results, which show that the proposed method significantly reduces order wait times compared to traditional methods. These findings underscore the advantages of integrating real-time adaptability and implicit forecasting into dynamic order-picking systems.

\section{Methodology} \label{sec:method}

\begin{figure}[!htbp]
\centering
\includegraphics[scale=0.45]{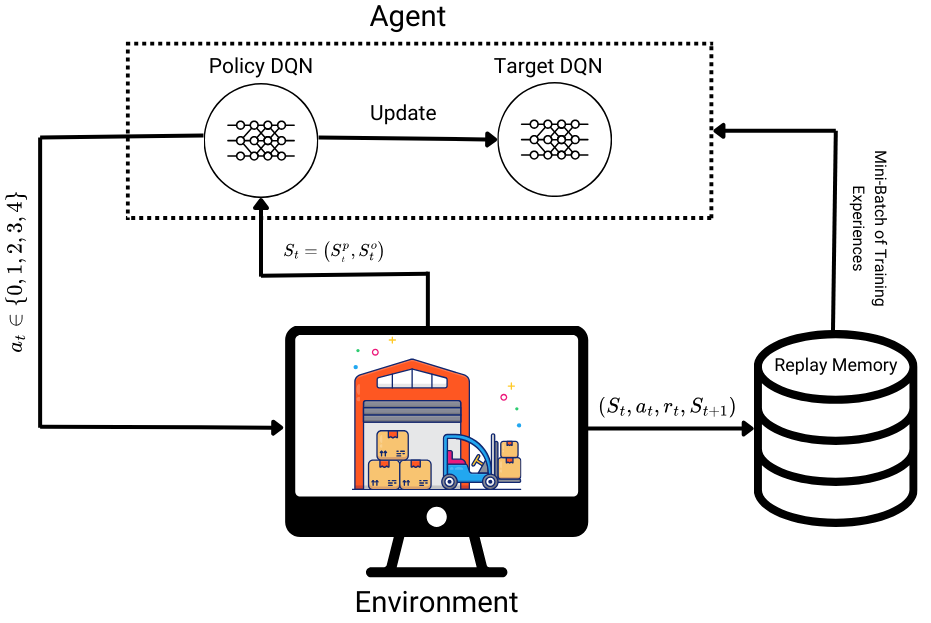}
\caption{A visual representation of the proposed solution method}
\label{fig:SummaryDQN}
\end{figure}

In this section, we present the formulation of the real-time warehouse order picking problem. We model this problem using a finite Markov Decision Process (MDP) with discrete time steps to capture its dynamic nature. Specifically, each time step $t$ involves making a movement or a decision regarding the picker's actions. At each step, the system state $S_t$ is observed, encapsulating information about both the orders, their locations in warehouse and the picker. Based on this state, an action $a_t$ is chosen, determining the picker's movement direction, whether the picker should remain stationary, or if items should be released at the depot. Executing this action leads to a new system state $S_{t+1}$, from which the next action $a_{t+1}$ is determined for the subsequent time step $t+1$. To solve the MDP model, we employ a custom-built Deep Q-learning approach. A high-level overview of this approach is illustrated in Figure~\ref{fig:SummaryDQN}, with detailed explanations provided in this section.

\subsection{MDP Formulation}
\label{subsec:MDP}
MDPs provide a structured framework for modeling decision-making processes where outcomes are influenced by both random factors and the actions of a decision maker. An MDP is defined by a five-tuple $(\textit{\textbf{S}}, \textit{\textbf{A}}, \textit{\textbf{P}}, \textit{\textbf{R}}, \gamma)$, where $\textit{\textbf{S}}$ represents the set of system states, 
$\textit{\textbf{A}}$ is a finite set of actions, $\textit{\textbf{P}}$ denotes the state transition probability, 
$\textit{\textbf{R}}$ is the immediate reward function, and $\gamma$ is the discount factor.  The following are the details of the specific formulation of the MDP for our problem.

\begin{enumerate}
    \item \textit{State}: The system state at time step $t$ is denoted as $\boldsymbol{S_t} = (\boldsymbol{S^p_t}, \boldsymbol{S^o_t})$, where $\boldsymbol{S^p_t}$ represents the picker state and $\boldsymbol{S^o_t}$ represents the state of existing orders in the warehouse. In the context of our research, using raw data to represent each state requires managing high-dimensional information. Specifically, representing the picker's state involves capturing its physical location coordinates on the grid and its current capacity. Similarly, representing existing orders from the raw data requires knowledge of the location coordinates for every storage location and the number of items ordered from each location, which could range from zero to a positive value. Consequently, the dimensionality of the state space derived from raw data can be exceedingly large. To address this complexity and facilitate decision-making in our research, we employ features to characterize the system's state. These features are designed based on two fundamental assumptions: When the vehicle's capacity is full, it must immediately return to the depot to drop off the items.
If the picker passes a pick location, it will collect all items available there, provided its capacity allows. Additionally, we assume that picking each item requires some time (depending on the application). With this framework in mind, we next present the features used to represent each state.
    
The picker state, $\boldsymbol{S^p_t} = (S^H_t, S^{V_1}_t, S^{V_2}_t, S^C_t)$, consists of four components: $S^H_t$, $S^{V_1}_t$, $S^{V_2}_t$, and $S^C_t$. Here, $S^H_t$ represents the horizontal position of the picker, which is the relative position of the picker with respect to cross-aisles, $S^{V_1}_t$ and $S^{V_2}_t$ represent the vertical position, which are picker's positions with respect to pick aisles, and $S^C_t$ indicates the remaining capacity of the picker. If $S^C_t = 0$, the picker is full and cannot pick up additional items; if $S^C_t = K$, the picker is empty. The horizontal state $S^H_t$ can take values of -1, 0, or 1: -1 indicates the picker is at the back cross-aisle, 1 indicates the front cross-aisle, and 0 indicates the picker is within one of the pick aisles. The exact location of the picker within a pick aisle when $S^H_t = 0$ is implicitly encoded in the order state $\boldsymbol{S^o_t}$. The vertical position is determined by $S^{V_1}_t$ and $S^{V_2}_t$, which equal $2n-1$ and $2n$ when the picker is at aisle $n$. 
For example, if the picker is inside aisle $4$ at time point $t$, then the vertical components of the picker state would be $S^{V_1}_t=7$ and $S^{V_2}_t=8$. The rationale for representing each pick aisle by two consecutive indices, rather than a single index, will be explained once the order state is described in more details. \\


The state of orders, $\boldsymbol{S^o_t} = (S^1_t, S^2_t, \ldots, S^{2N-1}_t, S^{2N}_t)$, where $N$ is the number of pick aisles in the warehouse. Each aisle $n$ is represented by $S^{2n-1}_t$ and $S^{2n}_t$. These values indicate the state of orders in each pick aisle depending on the picker's vertical movement direction. $S^{2n-1}_t$ represents the potential obtainable reward from aisle $n$ if the picker moves upward in its current aisle, while $S^{2n}_t$ represents the potential reward if the picker moves downward. These values are computed by summing the value of each order in the aisle divided by its distance from the picker, based on the direction of movement.  The distance between the picker and an order is the smallest number of movements required for picker to reach the order in terms of the number of location storage. 
These states are calculated as 
\begin{align}
    S^{2n-1}_t &= \displaystyle\sum_{i=1}^{L} \displaystyle \frac{n_i}{\rho^1_i} \\
    S^{2n}_t &=\displaystyle\sum_{i=1}^{L} \displaystyle \frac{n_i}{\rho^2_i},
\end{align}
where $L$ is the number of pick locations in each aisle, $n_i$ is the number of orders at pick location $i$, and $\rho^1_i$ and $\rho^2_i$ are the distances of picker to pick location $i$ when moving upward or downward, respectively.\\

For example, consider the picker's position and the requested items in Figure \ref{fig:warehouse}.
We assume that the picker is empty and thus the picker state is $\boldsymbol{S^p_t} = (0, 1, 2, $K$)$.
The distance traveled through each unit square in the figure is 1 unit, and the distance between adjacent cross-aisles is 3 units. 
Note that each storage location with a requested item in the figure contains only one requested item.
Now, the order state can be determined as follows.
For the first aisle, 
$$S^1_t = 0 \mbox{ and } S^2_t = \displaystyle\frac{1}{3}=0.33.$$ 
Similarly, for the second aisle, 
$$S_3=\displaystyle\frac{1}{9}+\displaystyle\frac{1}{15}=0.18 \mbox{ and } S_4=\displaystyle\frac{1}{13}+\displaystyle\frac{1}{19}=0.13.$$
Thus, the order state is calculated as $\boldsymbol{S^o_t} = (0, 0.33, 0.18, 0.13, 0, 0, 0.06, 0.04, 0.10, 0.08)$.\\

This formulation implicitly captures the importance of the picker's visiting direction for each aisle. Notably, when $S^H_t = -1$, then $S^{2n}_t = 0$ for all aisles since the picker cannot move upward. Similarly, when $S^H_t = 1$, then $S^{2n-1}_t = 0$ for all aisles. 
The starting state is $\boldsymbol{S^o_0} = \boldsymbol{0}$ since there are no orders at time zero, and $\boldsymbol{S^p_t} = (1, 2d-1, 2d, K)$, where $d$ is the aisle where depot is located at, and $K$ is the picker capacity. Finally, it is worth mentioning that since the order state $\boldsymbol{S^o_t}$ is formulated as a vector of size $2N$, we have intentionally chosen the same pattern for vertical components of the picker state for consistency.



    \item \textit{Action}: Given the state $S_t$, the action $a_t$ determines the picker's decision. The action $a_t$ can take on values from 0 to 4. Specifically, $a_t = 0$ signifies that the picker drops off items if located at the depot, and otherwise results in no movement. Actions with values 1 and 2 correspond to moving to the immediate next aisle, specifically to the right and left aisles, respectively. In this study, actions are constrained to feasible options only. For instance, when the picker is situated within a pick aisle (i.e., $S^H_t = 0$), actions 1 and 2 are disallowed. Actions 3 and 4 correspond to moving the picker upward and downward, respectively.  
The movement associated with Actions 3 and 4 is not necessarily limited to one unit; it should be sufficient to alter the state of the system as defined in the state transition explanation. Additionally, when the picker is in the back cross-aisle (i.e., $S^H_t = -1$), action 3 is prohibited, and when the picker is in the front cross-aisle (i.e., $S^H_t = 1$), action 4 is prohibited.

    \item \textit{State Transition}: In our MDP formulation, the state transition describes how the system evolves from one state to another based on the actions taken by the picker and the randomness involved in the order arrival. Let $\tau$ represent the smallest unit of time in which the picker can move from one storage location to another. If action $a_t = 0$ is executed at any location other than the depot, it indicates that the picker remains stationary for $\tau$ time until the next step. Consequently, $\boldsymbol{S^p_{t+1}}$ will be identical to $\boldsymbol{S^p_t}$, and $\boldsymbol{S^o_{t+1}}$ will remain unchanged unless a new order arrives during $\tau$, in which case it will be recalculated. The number of new orders arriving within this interval depends on the order arrival rate.\\ 
    
    If action $a_t = 0$ is executed at the depot, the picker will unload all items. Therefore, if $\boldsymbol{S^p_t} = (1, 2d-1, 2d, S^C_t)$ and $a_t = 0$, then $\boldsymbol{S^p_{t+1}} = (1, 2d-1, 2d, K)$ and $\boldsymbol{S^o_{t+1}}$ will be recalculated based on the orders arriving during the time required to unload $K - S^C_t$ items. Actions $a_t = 1$ and $a_t = 2$ can only be performed in cross aisles and require time proportional to the inter-aisle distance. Given $\boldsymbol{S^p_t} = (S^H_t, S^{V_1}_t, S^{V_2}_t, S^C_t)$, if $a_t = 1$ is executed, the new state will be $\boldsymbol{S^p_{t+1}} = (S^H_t, S^{V_1}_t + 2, S^{V_2}_t + 2, S^C_t)$.
Conversely, if $a_t = 2$ is executed, $\boldsymbol{S^p_{t+1}} = (S^H_t, S^{V_1}_t - 2, S^{V_2}_t - 2, S^C_t)$.\\

In our MDP formulation, actions 3 and 4 result in a new state under one of the following conditions: the picker reaches a cross aisle and cannot proceed further, a new order arrives altering $\boldsymbol{S^o_t}$, or the picker collects at least one item. It is assumed that the picker collects all items at a storage location before transitioning to a new state. 
When action $a_t = 3$ or $a_t = 4$ is executed, the state of orders $\boldsymbol{S^o_t}$ will change irrespective of new order arrivals due to the altered relative distance caused by the picker's movement. For the picker state $\boldsymbol{S^p_{t+1}}$, $S^H_{t+1}$ may change based on whether the picker moves from a cross aisle to within an aisle or vice versa, but $S^{V_1}_{t+1}$ and $S^{V_2}_{t+1}$ will remain unchanged. The value of $S^C_{t+1}$ will depend on the number of items collected during the action.

    \item Reward: The reward at time step $t$ is defined as follows:

\begin{equation}
r_t = 
\begin{cases} 
-1 & \text{if } a_t = 0 \text{ and picker not at the depot}, \\
R \times (K - S^C_t) \times \alpha & \text{if } a_t = 0 \text{ and picker at the depot}, \\
- n^m_t + R \times n^p_t & \text{otherwise}.
\end{cases}
\end{equation}

\textcolor{black}{
The components of the reward function are interpreted as follows:
\begin{enumerate}
    \item Case: \( a_t = 0 \) and the picker is not at the depot.
This action signifies that the picker remains stationary at its current location. A negative reward of \(-1\) is assigned in this case, incentivizing the picker to move when a potential positive reward can be obtained rather than staying idle. This penalty aligns with the penalty for moving one storage location, thereby encouraging exploration without overly penalizing idle actions. By making the penalty for staying still equal to the cost of moving, the policy is encouraged to weigh the potential rewards of movement more effectively.
\item Case: \( a_t = 0 \) and the picker is at the depot.
This action represents unloading all items from the picker. The reward for this action is proportional to the number of items being unloaded, which is given by \( K - S^C_t \), where \( S^C_t \) denotes the current capacity of the picker. The parameter \( \alpha \) modulates the reward for unloading items relative to picking them up, with its value ranging between 0 and 1. Smaller values of \( \alpha \) prioritize picking up items over unloading, encouraging the picker to maximize its capacity before heading to the depot, which helps reduce travel distance. Conversely, larger values of \( \alpha \) balance the value of unloading and picking, which minimizes both waiting and delivery times.
For example, suppose \( R = 8 \), \( K = 10 \), \( S^C_t = 7 \), and \( \alpha = 0.5 \). In this case, the picker unloads \( 3 \) items at the depot, receiving a reward of:
\[
R \times (K - S^C_t) \times \alpha = 8 \times 3 \times 0.5 = 12.
\]
\item Case: \( a_t \neq 0 \).
Any action other than \( 0 \) leads to the picker moving to a new location. The reward in this case depends on both the distance traveled (number of storage locations) and the number of items picked up during this movement.
\begin{itemize}
    \item The number of storage locations traversed is penalized by a term \( n^m_t \), representing the effort or cost of movement.
    \item If items are picked up along the way, the reward is augmented by a term \( R \times n^p_t \), where \( n^p_t \) is the number of items picked up.
    \item For actions \( a_t = 1 \) and \( a_t = 2 \), which correspond to moving left and right within the cross aisles, no items are picked up. Hence, the reward includes only the penalty for movement (\(- n^m_t\)).
    \item For actions \( a_t = 3 \) and \( a_t = 4 \), which involve moving and picking up orders, the reward reflects both penalties for movement and the reward for picking items.  
    For example, suppose the picker moves 7 storage locations and picks up 2 items, with \( R = 8 \). The reward is calculated as:
    \[
    - n^m_t + R \times n^p_t = -7 + 8 \times 2 = 9.
    \]
\end{itemize}
\end{enumerate}
}

\end{enumerate}

Additional discussion regarding the rationale for the feature-based state representation proposed in this study is provided in \ref{app:rationale}.

\begin{figure}
    \centering
    \includegraphics[width=\textwidth]{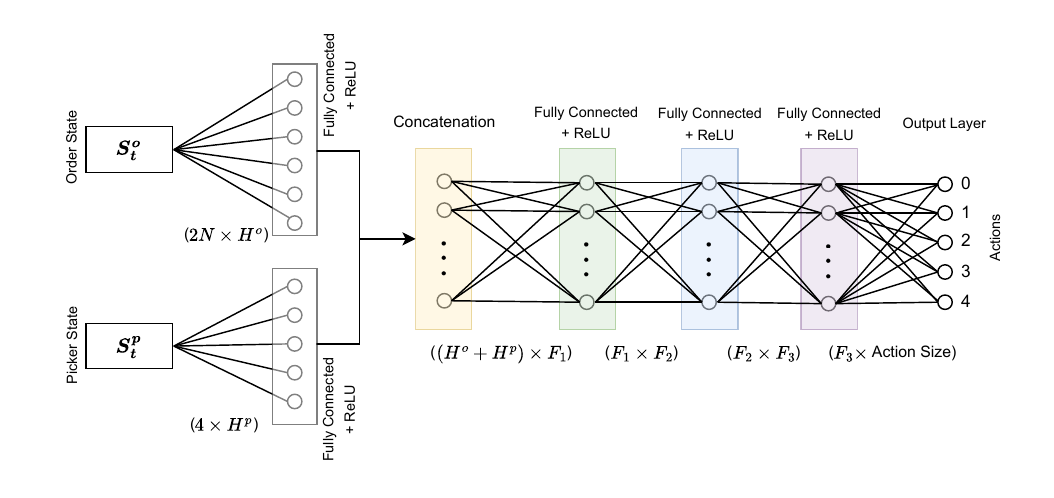}
    \caption{Proposed Deep Neural Network Architecture}
    \label{fig:DNN}
\end{figure}

\subsection{Proposed Approach}
Our objective is to train a policy that aims to maximize the discounted cumulative reward $R_t = \sum_{k=0}^{\infty} \gamma^k \times r_{t+k}$, where $R_t$ is referred to as the return, and $0 < \gamma < 1$ is the discount factor that balances the importance between immediate and future rewards. Small values of $\gamma$ indicate that immediate actions are more significant than future actions, while larger values of $\gamma$ place greater emphasis on future rewards.

The fundamental concept behind Q-learning \citep{watkins1989learning} is that if we had access to a function $Q^*: \text{State} \times \text{Action} \rightarrow \mathbb{R}$, which could inform us of the expected return if we were to take a particular action in a given state, we could straightforwardly construct a policy that maximizes our rewards: $\pi^* = \arg\max_a Q^*(s, a)$, where $Q^*(s, a) = \max (E_{\pi^*}[R_t \mid s_t = s, a_t = a])$.  However, due to the complexity of the environment, we do not have direct access to $Q^*$. Nevertheless, since neural networks are universal function approximators, we can construct a neural network and train it to approximate $Q^*$.

\subsubsection{Architecture of Deep Neural Network} 

The deep neural network architecture designed for this study is specifically tailored to handle the state representation and predict optimal actions for the picker in a dynamic warehouse environment. The network serves as a Q-network, approximating the action-value function \( Q(s, a) \), and consists of several fully connected layers to process both the picker state and the order state separately before combining them for the final action prediction. Figure \ref{fig:DNN} provides a summary of the architecture and it is structured as follows:
\begin{itemize}
\item \textit{Input Layers:}
       The input to the neural network is divided into two parts: the picker state \( S^p_t \) and the order state \( S^o_t \). The picker state \( S^p_t \) is a four-dimensional vector represented as \( (S^H_t, S^{V_1}_t, S^{V_2}_t, S^C_t) \), while the order state \( S^o_t \) is a vector of size \( 2N \), where \( N \) is the number of aisles in the warehouse. \textcolor{black}{Note that the input to the neural network is the state representation as defined earlier, comprising both picker state ($S^p_t$) and order state ($S^o_t$). This ensures that the neural network processes the exact state information defined for the problem.}

\item \textit{Processing Picker State:}
 The picker state input \( S^p_t \) is processed by the first fully connected layer, transforming the four-dimensional picker state vector into a \( H^p \)-dimensional feature vector using a ReLU activation function. This transformation helps in capturing the essential features of the picker state.

\item \textit{Processing Order State:}
The order state input \( S^o_t \) is processed by a second fully connected layer, which transforms the \( 2N \)-dimensional order state vector into a \( H^o \)-dimensional feature vector using a ReLU activation function. This layer extracts meaningful features from the state of orders in the warehouse.

\item \textit{Combining Picker and Order States:}
 The feature vectors obtained from processing the picker state and order state are concatenated into a single vector. This combined vector with size $H^p+H^o$, integrates the information from both states, resulting in a comprehensive representation of the system state.

\item \textit{Further Processing of Combined Features:}
 The concatenated vector is then passed through a series of fully connected layers:
       \begin{itemize}
           \item The first layer in this series has \( F_1 \) units and uses a ReLU activation function to further transform the combined feature vector.
           \item The subsequent layer reduces the dimensionality to \( F_2 \) units, again using a ReLU activation function.
           \item Another layer reduces the dimensionality further to \( F_3 \) units, continuing with the ReLU activation function.
       \end{itemize}

\item \textit{Output Layer:}
 Finally, the output layer maps the \( F_3 \)-dimensional feature vector to the action space, predicting the Q-values for each possible action. The action with the highest Q-value is then selected as the optimal action for the picker.

\end{itemize}
\textcolor{black}{We utilize this deep neural network architecture to capture the complex interactions between the picker and the orders in the warehouse, enabling it to predict optimal actions that maximize the cumulative reward over time.}

\subsubsection{Comparison with Alternative Neural Network Architectures}
\label{subsubsec:nn_comparison}


To illustrate why a fully connected feed-forward neural network (FCN) is particularly well-suited to our problem setting, we compare it with other common neural network architectures. Table~\ref{tab:arch_comparison} summarizes key architectural choices, their typical computational complexity, and their applicability to the warehouse order-picking domain.

\begin{table}[ht]
\scriptsize
\centering
\caption{Comparison of common neural network architectures}
\label{tab:arch_comparison}
\begin{tabular}{lll}
\toprule
\textbf{Architecture} &
\begin{tabular}[l]{@{}l@{}}\textbf{Computational Complexity} \end{tabular} &
\begin{tabular}[l]{@{}l@{}}\textbf{Suitability}\end{tabular}  \\
\midrule
\begin{tabular}[l]{@{}l@{}}\textbf{Fully Connected} \\ \textbf{Network (FCN)}\\ (Our Approach)\end{tabular} 
 & 
\begin{tabular}[l]{@{}l@{}}$\mathcal{O}(n \cdot m)$ per layer, \\ where $n$ = \# input features, \\ $m$ = \# outputs.\end{tabular}
&
\begin{tabular}[l]{@{}l@{}}Well-suited for structured, \\ tabular-like input\end{tabular} \\
\hline
\begin{tabular}[l]{@{}l@{}}\textbf{Convolutional Neural} \\ \textbf{Network (CNN)}\end{tabular}
 &
\begin{tabular}[l]{@{}l@{}}$\mathcal{O}(H \cdot W)$ per layer, \\ where $H \times W$ = 2D input size. \\ Assumes fixed kernel and \\ channel sizes.\end{tabular}
&
\begin{tabular}[l]{@{}l@{}}Ideal for spatial/\\image-based data\end{tabular} \\
\hline
\begin{tabular}[l]{@{}l@{}}\textbf{Graph Neural} \\ \textbf{Network (GNN)}\end{tabular}
&
\begin{tabular}[l]{@{}l@{}}$\mathcal{O}(|V| + |E|)$ per layer, \\ where $|V|$ = \# nodes, $|E|$ = \# edges.\\ Assumes sparse graph and fixed \\ feature dimensions.\end{tabular}
&
\begin{tabular}[l]{@{}l@{}}Useful for explicit graph\\representations\end{tabular} \\
\hline
\begin{tabular}[l]{@{}l@{}}\textbf{Recurrent Neural} \\ \textbf{Network (RNN)} \end{tabular}
&
\begin{tabular}[l]{@{}l@{}}$\mathcal{O}(T)$ per layer, \\ where $T$ = sequence length. \\ Assumes fixed hidden and input \\ dimensions across steps.\end{tabular}
&
\begin{tabular}[l]{@{}l@{}}Effective for temporal\\or sequential data\end{tabular} \\
\bottomrule
\end{tabular} 
\end{table}

Given the aggregated vector-based state representation introduced in Section~\ref{subsec:MDP}, the proposed framework utilizes two types of inputs: a relatively low-dimensional picker state vector \((S^H_t, S^{V_1}_t, S^{V_2}_t, S^C_t)\) and a higher-dimensional, yet still tabular, order state vector \((S^1_t, S^2_t, \ldots, S^{2N}_t)\), which captures aggregated aisle-level demand. These inputs constitute structured tabular data, as opposed to raw image-based inputs or explicit spatial grids. Accordingly, a fully connected network (FCN) is adopted, given its effectiveness in processing tabular data with relatively low computational overhead. This makes the architecture well-suited for dynamic and time-sensitive tasks such as real-time order-picking decisions. Moreover, spatial relationships and capacity constraints are already accounted for through targeted feature engineering---such as aisle-level demand aggregation and distance-based reward shaping---thereby mitigating the need for more specialized architectures like convolutional neural networks (CNNs), graph neural networks (GNNs), or recurrent neural networks (RNNs). These alternatives are typically employed when handling data with strong spatial, graphical, or temporal dependencies, which are not explicitly present in our setting.

\color{black}

\subsubsection{Training of Deep Neural Network}

The training process leverages two neural networks: a policy network and a target network. 
The policy network guides decision-making, selecting actions based on the epsilon-greedy strategy, which balances exploration of action space and exploitation of the learned values.
The target network, on the other hand, serves as a stable reference, updated periodically with the policy network's learned parameters.
To enhance the learning process, we employ the replay memory mechanism.
Replay memory acts as a repository to store the agent's experiences in the form of state transitions. 
During training, the DRL agent randomly samples mini-batches of transitions from the replay memory. 
This allows for efficient reuse of past experiences, which helps to stabilize the learning process. 

\begin{algorithm}[!htbp] 
\caption{Training Procedure for DQN in Real-Time Warehouse Order Picking}
\label{alg:DQN}
\begin{algorithmic}[1] 
    \State Initialize policy network $\pi_{\theta}$ and target network $\pi_{\theta^-}$ with weights $\theta$ and $\theta^- = \theta$
    \State Initialize replay memory $D$ with capacity $C$
    \For{episode = 1 to $N_e$}
        \State Initialize state $S_0$
        \For{t = 1 to $N_s$}
            \State With probability $\epsilon$, select a random action $a_t$
            \State Otherwise, select $a_t = \arg\max_a Q(S_t, a; \theta)$
            \State Execute action $a_t$ and observe reward $r_t$ and next state $S_{t+1}$
            \State Store transition $(S_t, a_t, r_t, S_{t+1})$ in replay memory $D$
            \If{size of $D$ $\geq$ mini-batch size}
                \State Sample random mini-batch of transitions from $D$
                \State Set $S_M$, $S_M'$ as states and next states in the mini-batch
                \State Set $A$ as actions, $R$ as rewards
                \State Compute $Q(S_M, A; \theta)$
                \State Set $Q_{\text{target}} = R + \gamma \times \max_{a' } Q(S_M', a'; \theta^-)$
                \State Compute loss $L = \text{HuberLoss}(Q(S_M, A; \theta), Q_{\text{target}})$
                \State Perform gradient descent on $L$ to update $\theta$
            \EndIf
            \If{$t \pmod{N_{\text{update}}} == 0$}
                \State Update target network $\theta^- = \tau \theta + (1 - \tau) \theta^-$
            \EndIf
        \EndFor
    \EndFor
\end{algorithmic}
\end{algorithm}

The training is conducted over \(N_e\) episodes, each comprising \(N_s\) steps. 
At each step, the agent observes the current state, selects an action using the epsilon-greedy method and the policy network, and and receives a reward upon transitioning to the next state. 
The current state, action, next state, and reward form a transition, and is stored in the replay memory with capacity $C$.
Subsequently, a randomly selected mini-batch of transitions is sampled from the memory, and the policy network undergoes one step of training based on these experiences. 
If the replay memory lacks sufficient transitions to form a mini-batch, the training step is skipped temporarily.

The training procedure for each mini-batch is a multi-step process. 
First, the states, next states, actions, and rewards within the batch are stacked for parallelization. 
Then, the policy network is set to training mode, while the target network enters evaluation mode. 
This setup allows for the computation of state-action values using the policy network, and expected state-action values using the target network.
The discrepancy between these values is then used to calculate the loss, and we utilize the Huber loss function \citep{huber1992robust} for robustness. 
The calculated loss is then backpropagated through the policy network, facilitating its optimization.
To ensure stability, the target network is updated less frequently, specifically every \(N_{\text{update}}\) steps.
This update follows the soft update rule using target update weight $\tau$, where the new target network weights are updated as the weighted average of the policy network weights and the existing target network weights.
This gradual update mechanism helps prevent drastic fluctuations and promotes smoother learning.
The Deep Q-Network (DQN) training procedure is outlined in the pseudocode presented in Algorithm~\ref{alg:DQN}.

\section{Computational Study} \label{sec:comp_exp}

\paragraph{Experimental Settings} The warehouse layout and instance generation process employed in this study are derived from \citet{lu2016algorithm}. 
The warehouse configuration consists of a single block with two cross-aisles and 10 parallel aisles, each containing 15 storage locations on either side. 
Each aisle is 15 meters long, with 3-meter spacing between adjacent aisles.
We assume that the depot is located at the end of aisle 6, adjacent to the front cross-aisle, with negligible travel time between these points.
The order-picker operates at a speed of 1 meter per second, with a picking time of 5 seconds per item and an additional 1 second for an item drop-off at the depot. 
The picking device has a capacity ($K$) of 20 items. 
We assume that orders arrive following a Poisson process and are distributed uniformly throughout the warehouse. We conduct experiments with arrival rates ranging from 0.01 to 0.09 orders per second. 
For each order arrival rate, we conduct 10 independent simulation runs and report the average results to ensure reliable performance assessment.
Each experiment simulates an 8-hour work shift, during which incoming orders are dynamically assigned to the picker to optimize warehouse order-picking throughput. 
An order is considered complete only after the picker deposits the item at the depot, and not upon retrieval from the respective storage location.

\paragraph{Evaluation Metrics} To evaluate the efficiency of the proposed dynamic order-picking algorithm, we employ three key metrics:
\begin{enumerate}
    \item \textit{Average Travel Distance per Order (ATDO):} 
    The total distance traveled by an order picker during a work shift, divided by the number of completed orders.
    \item \textit{Average Order Completion Time (AOCT):} 
    The average time elapsed between an order arrival at the warehouse and its final delivery at the depot by the order-picker. Our research prioritizes minimizing AOCT. By streamlining order processing, we aim to reduce customer wait times and consequently improve order fulfillment rates.
    \item \textit{Percentage of Unfulfilled Orders (PUO):} 
    The percentage of orders that are not successfully completed within the work shift. A high PUO may indicate high order arrival rates, limited picking device capacity, or inefficient routing. 
\end{enumerate}
Traditionally, when designing order picking algorithms or strategies, more focus has been placed on the shortest path. 
However, as highlighted in Section~\ref{sec:lit_rev}, the primary objective in warehouse management is to minimize order throughput time. 
Therefore, our solution methodology prioritizes both AOCT and PUO.


\subsection{Benchmark Algorithms}
To evaluate our order-picker routing algorithm, we employ a total of seven benchmark approaches. Two of these are widely used heuristic picking policies from the literature, specifically the S-shape and Largest-gap methods, which we discussed earlier in Section~\ref{sec:lit_rev}. These serve to analyze the solution quality of our algorithm on the problem at hand. The remaining five benchmarks are based on different configurations of approaches introduced in two well-known studies from the literature:
\begin{enumerate}
    \item \citet{ratliff1983order}: This seminal work introduced an efficient algorithm for the static order picking problem.
    The order picker operates under a batching strategy, waiting for the pick list to reach the capacity of picking device, followed by a shortest tour to fulfill the orders. 
    The authors employ a graph-based dynamic programming approach to reduce the computational complexity compared to standard TSP solution methods.
    \item \citet{lu2016algorithm}: This study adapts the static order picking approach for dynamic environments, under certain assumptions. 
    Specifically, the order picker begins a tour when incoming orders reach 25\% of device capacity and can be re-routed if new orders arrive within capacity limits. 
    Furthermore, assuming narrow cross-aisles, the picker cannot be re-routed while traveling within them. 
    This method is referred to as the Interventionist Routing Algorithm (IRA).
    The authors have demonstrated the superiority of their method compared to an effective heuristic known as the largest-gap heuristic.
\end{enumerate} 

We note that the source code for the aforementioned algorithms is unavailable. More importantly, the core of the above-mentioned studies relies on Dynamic Programming (DP) methods, which are favored for their efficacy in picker routing. However, these DP methods operate under certain assumptions, and extending them to accommodate the diverse routing scenarios found in real-world warehouse operations is not a trivial task. For instance, although picker re-routing within cross-aisles is practical and significantly impacts operational efficiency (as will be shown in the numerical results), it is assumed to be prohibited in these DP methods. Note that our proposed DRL agent automatically learns to apply cross-aisle re-routing as needed. Therefore, to replicate the results of the aforementioned studies under their assumptions, as well as to make these methods adaptable for various routing scenarios, we replace DP with a generic TSP solver in their methods to compute optimal routes. This allows for more flexible routing scenarios, such as using smaller initial pick list sizes or permitting cross-aisle re-routing. 
Although a generic TSP solver is slower, it does not negatively impact the algorithm's decision-making process. This is because we exclude the computational time of the generic TSP solver in our experiments. In other words, during the simulation of work shifts, we assume the solution time of these generic solvers is zero and do not include it in the decision-making process.
Additionally, in our computational study, we do not report any computational time, as the proposed DRL method makes decisions instantaneously at any given state. Thus, the runtime is negligible during testing, with all computational time spent during the training phase. With these considerations, the remaining five baseline models for our computational study, based on optimal TSP routing, can be expressed as follows:

\begin{enumerate}
    \item \textbf{Baseline 1}: The initial pick list size $k$ is set equal to the picker capacity. This replicates the static routing approach of \citet{ratliff1983order}. Note that cross-aisle re-routing is not applicable in this scenario, as the picker's load capacity is already fully assigned when a new order arrives.
    \item \textbf{Baseline 2}: In this model, $k$ is set to 5, i.e., 25\% of the picker capacity, and the picker routing is performed without cross-aisle re-routing. This model represents the IRA in \citet{lu2016algorithm}.
    \item \textbf{Baseline 3}: A variant of the IRA where cross-aisle re-routing is enabled.
    \item \textbf{Baseline 4}: The initial pick list size is set to 1, meaning the picker departs from the depot as soon as the first order arrives. Cross-aisle re-routing is not permitted in this model.
    \item \textbf{Baseline 5}: Similar to Baseline 4, we set $k=1$. However, cross-aisle re-routing is enabled in this model.
\end{enumerate}

These variants serve as additional benchmarks to assess the effects of initial pick list size and cross-aisle re-routing constraints on overall picking performance.
Note that all baseline comparisons utilize a first-come-first-served batching policy, consistent with the above referenced benchmark studies. 
A comparison of different baseline assumptions is listed in Table \ref{tab:baselines}. We note that in this table, we report only one configuration of our proposed method, as only this configuration's results will be presented throughout this section. However, it is important to note that all the configurations used for Baselines 1–5 can also be applied to our proposed method. For interested readers, a comprehensive analysis of different variations of our proposed method is provided in \ref{app:analysis}.

\begin{table}[!htbp]
\scriptsize
    \centering
    \caption{Comparison of Model Assumptions}
    \label{tab:baselines}
    \begin{tabular}{ccc}
    \toprule
    Model      & Initial Pick Size  & Cross-aisle Re-routing \\
    \midrule
    Baseline 1 & 20 & N/A                    \\
    Baseline 2 & 5  & No                     \\
    Baseline 3 & 5  & Yes                    \\
    Baseline 4 & 1  & No                     \\
    Baseline 5 & 1  & Yes                    \\
    \color{black}{
    S-shape} & \color{black}{1}  & \color{black}{Yes}                     \\
    \color{black}{Largest-gap} & \color{black}{1}   & \color{black}{Yes}                     \\
    \hline
    Proposed DRL Agent  & 1  & Yes                   \\
    \bottomrule
    \end{tabular}
\end{table}

\subsection{Training and Testing Details}

\paragraph{Training} The proposed DRL approach is implemented in Python 3.9 using PyTorch 1.11. Interested readers can find our codes available on \href{https://github.com/Sasanm88/DRL-Order-Picking}{https://github.com/Sasanm88/DRL-Order-Picking}.  
The hyperparameter selection was conducted using a grid search, where the optimal value was determined based on the reward function. Table~\ref{tab:params} presents the results of the grid search for hyperparameters.
The selected configuration consists of a batch size of 64, with training conducted over 4,500 episodes, each comprising 1,000 steps. 
The neural network architecture consists of initial layers with dimensions $H^o = 160$ and  $H^p = 64$, followed by the fully connected layers with dimensions $F_1 = 256, F_2 = 128$ and $F_3 = 64$.
Note that $H^o$ and $H^p$ represent two components of a single layer, meaning that the structure $[H^o , H^p, F_1, F_2, F_3]$ effectively constitutes four layers in total.
Experience replay is employed with a memory capacity of $C = 200,000$. 
The optimization process utilizes a learning rate of 0.0001, discount factor $\gamma = 0.99$, and target update coefficient $\tau = 0.001$. 
The reward function, $R$, is defined as the sum of the number of storage locations in each aisle and the total number of aisles, ensuring a positive reward for successful item pickup regardless of the agent's initial location.

\begin{table}[!htbp]
    \centering
    \scriptsize
    \caption{Selected Training Hyperparameters Using Grid Search for the DRL Model}
    \label{tab:params}
    \begin{tabular}{lll}
    \toprule
    \textbf{Hyperparameter}      & \textbf{Grid}  & \textbf{Selected Value} \\
    \midrule
    No. of hidden layers & [1,2,3,4,5] & 4                   \\
    \midrule
    Hidden layer sizes & [[64, 32, 128, 256, 128], [128, 128, 128, 64, 32], & [160, 64, 256, 128, 64]                   \\
    $[H^o , H^p, F_1, F_2, F_3]$ & [128, 64, 64, 256, 128], [160, 64, 256, 128, 64]] &                      \\
    \midrule
    Activation & [ReLU, Tanh]  & ReLU                    \\
    \midrule
    Output activation & [Sigmoid, Softmax, Linear]  & Softmax                     \\
    \midrule
    Optimizer & [Adam, SGD]  & Adam                    \\
    \midrule 
    Batch size & [32, 64, 128, 256] & 64 \\
    \bottomrule
    \end{tabular}
\end{table}

\begin{figure}[ht]
    \begin{subfigure}[b]{0.49\textwidth}
        \centering
        \includegraphics[width=\textwidth]{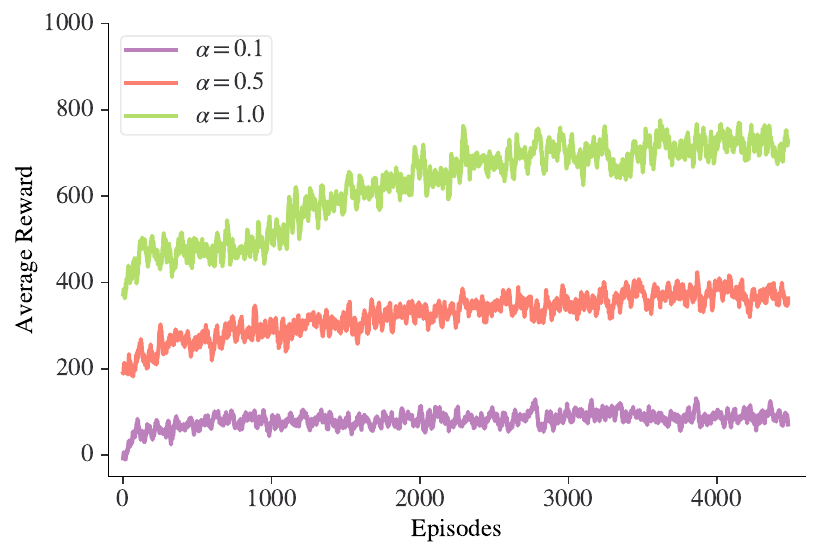}
        \caption{$\lambda = 0.02$}
        \label{subfig:lambda2}
    \end{subfigure}
    \hfill
    \begin{subfigure}[b]{0.49\textwidth}
        \centering
        \includegraphics[width=\textwidth]{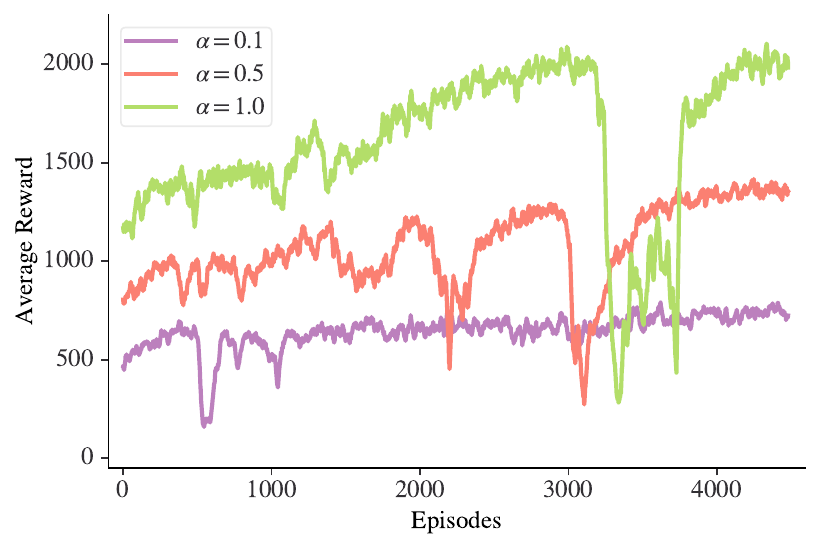}
        \caption{$\lambda = 0.04$}
        \label{subfig:lambda4}
    \end{subfigure}
    \\
    \begin{subfigure}[b]{0.49\textwidth}
        \centering
        \includegraphics[width=\textwidth]{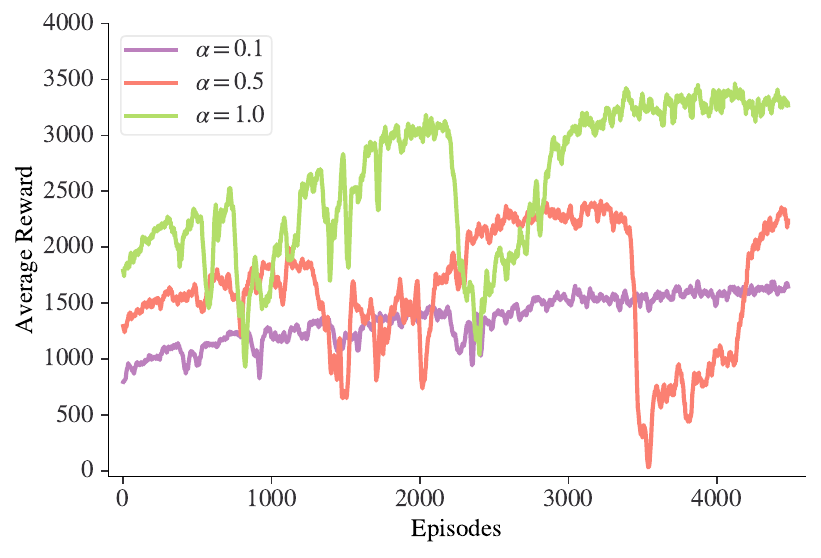}
        \caption{$\lambda = 0.06$}
        \label{subfig:lambda6}
    \end{subfigure}
    \hfill
    \begin{subfigure}[b]{0.49\textwidth}
        \centering
        \includegraphics[width=\textwidth]{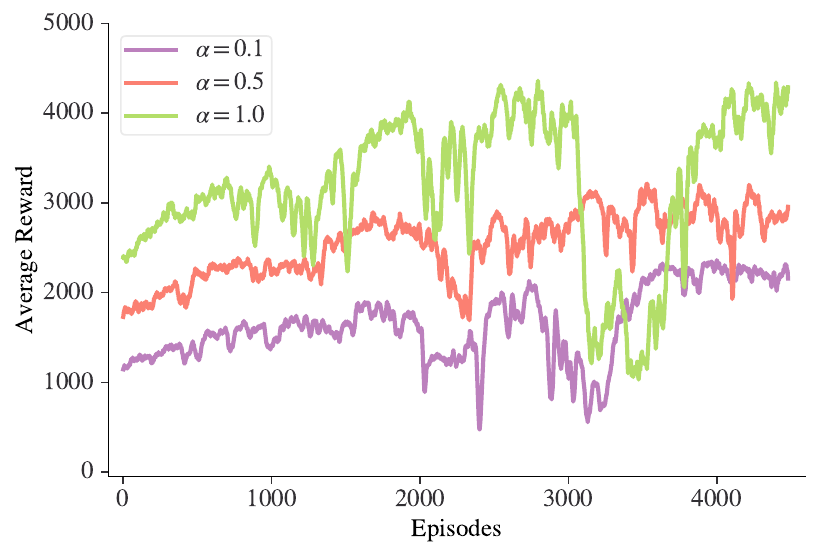}
        \caption{$\lambda = 0.08$}
        \label{subfig:lambda8}
    \end{subfigure}
    \caption{Training curves showing the convergence of the 20-period moving average reward for each trained model}
    \label{fig:training_curves}
\end{figure}

We note that a DRL model trained for a specific order arrival rate may not generalize well to different arrival rates. Machine learning frameworks often exhibit sensitivity to the training distribution and its parameter values \citep{tsotsos2019does, sobhanan2023genetic}. So, we conduct a series of training experiments with varying arrival rates: $\lambda \in \{0.02, 0.04, 0.06, 0.08\}$.
Furthermore, we investigate the trade-off between prioritizing distance minimization and order completion time by adjusting the weight hyperparameter $\alpha \in \{0.1, 0.5, 1.0\}$. This results in the training of 12 models, given that there are four values for $\lambda$ and three values for $\alpha$. 
The training process for each model takes approximately 7 hours to complete on a MacBook Pro equipped with an Apple M1 Chip and 16GB RAM. 
The training curves for each trained model with varying $\lambda$ and $\alpha$ are shown in Figure~\ref{fig:training_curves}. These figures illustrate the convergence of the training procedure. The training curves indicate that the agent tends to engage in exploration when there are no significant changes in the average reward over multiple episodes.

\paragraph{Testing and Generalibility} We evaluate the performance of each trained model using test instances with varying values of $\lambda$ to assess their generalizability across different arrival patterns and optimization objectives. For this testing, we use 9 distinct values for $\lambda$, ranging from $0.01$ to $0.09$ in increments of $0.01$, and evaluate each of the 12 trained models against all of these values. Each model is tested against 1 in-sample scenario and 8 out-of-sample scenarios. Note that test scenarios with $\lambda \in \{0.01, 0.03, 0.05, 0.07, 0.09\}$ serve as out-of-sample scenarios for all trained models. As previously mentioned, each testing scenario consists of 10 simulation runs, with each run representing an 8-hour work shift. It is worth mentioning that in addition to the experiments presented in the main body of the paper, we also conducted supplementary simulations, presented in \ref{app:analysis_varying_rates}, where the arrival rate varies over the course of the 8-hour shift. The performance of the proposed algorithms under these time-varying arrival patterns was consistent with the results reported in the main analysis, where a constant arrival rate was assumed throughout the shift. These findings further support the robustness of the proposed models.
  
\begin{figure}[ht]
    \begin{subfigure}[b]{0.49\textwidth}
        \centering
        \includegraphics[width=\textwidth]{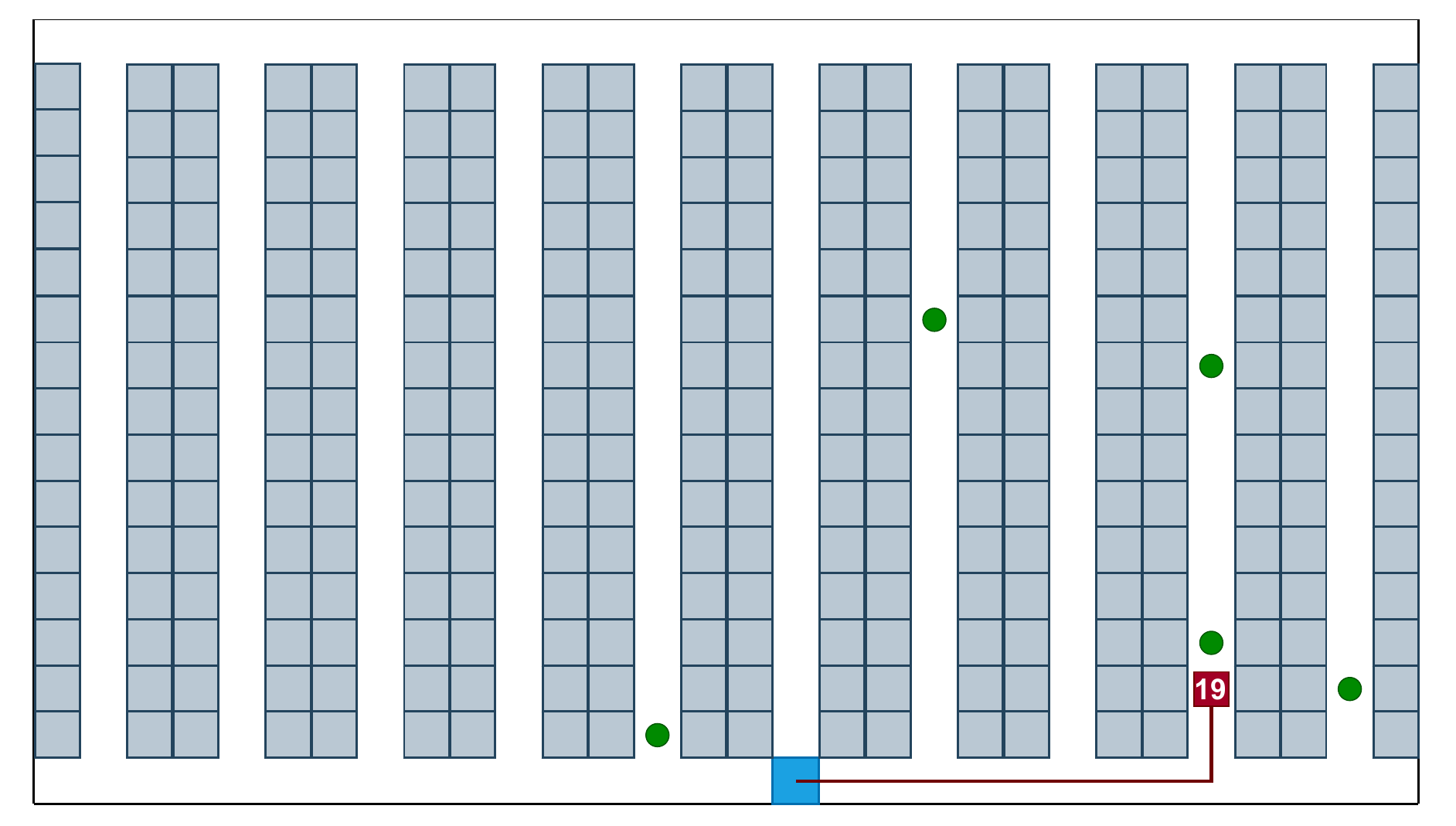}
        \caption{Time (sec.): 673}
        \label{subfig:673}
    \end{subfigure}
    \hfill
    \begin{subfigure}[b]{0.49\textwidth}
        \centering
        \includegraphics[width=\textwidth]{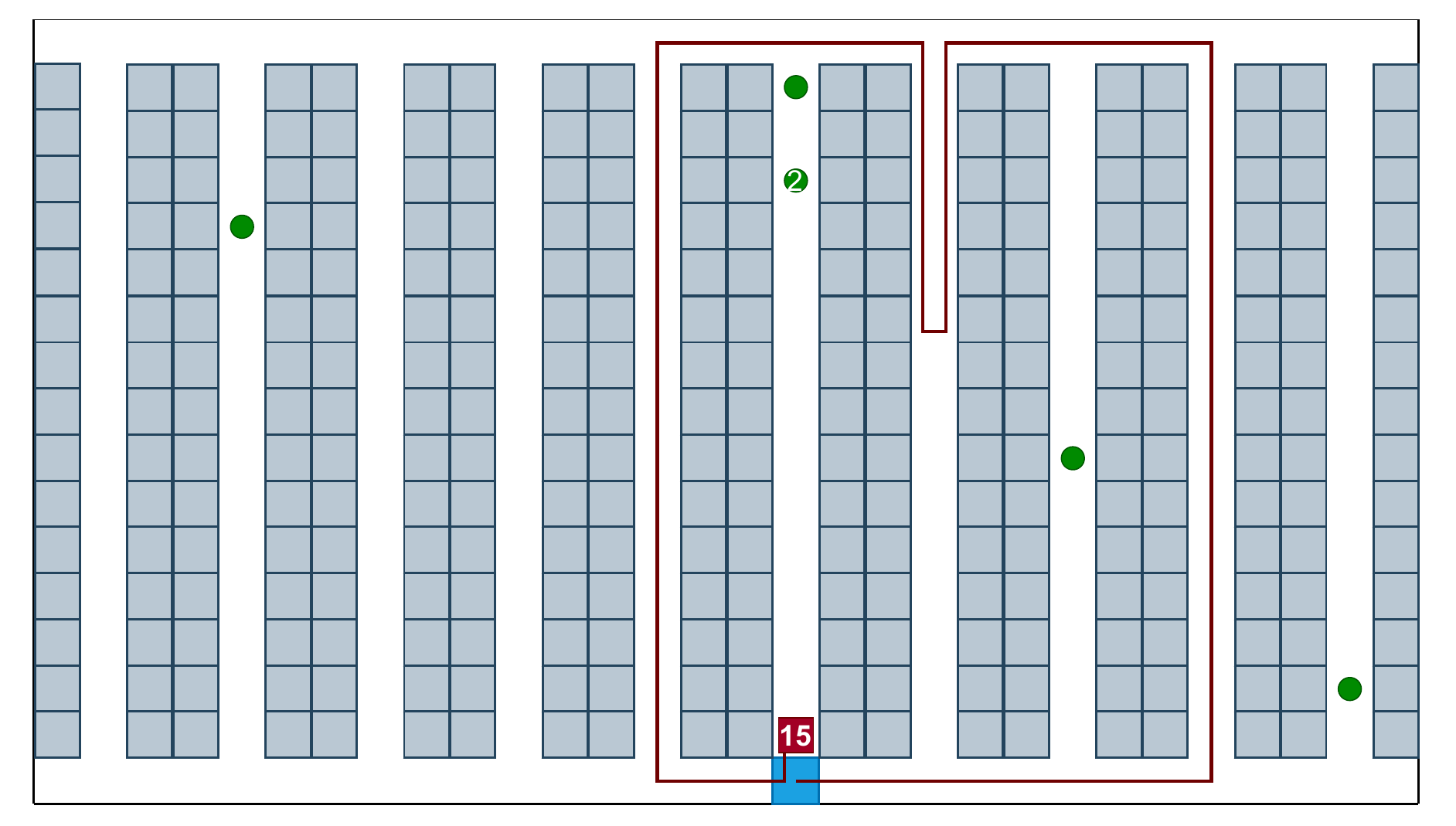}
        \caption{Time (sec.): 735}
        \label{subfig:735}
    \end{subfigure}
    \\
    \begin{subfigure}[b]{0.49\textwidth}
        \centering
        \includegraphics[width=\textwidth]{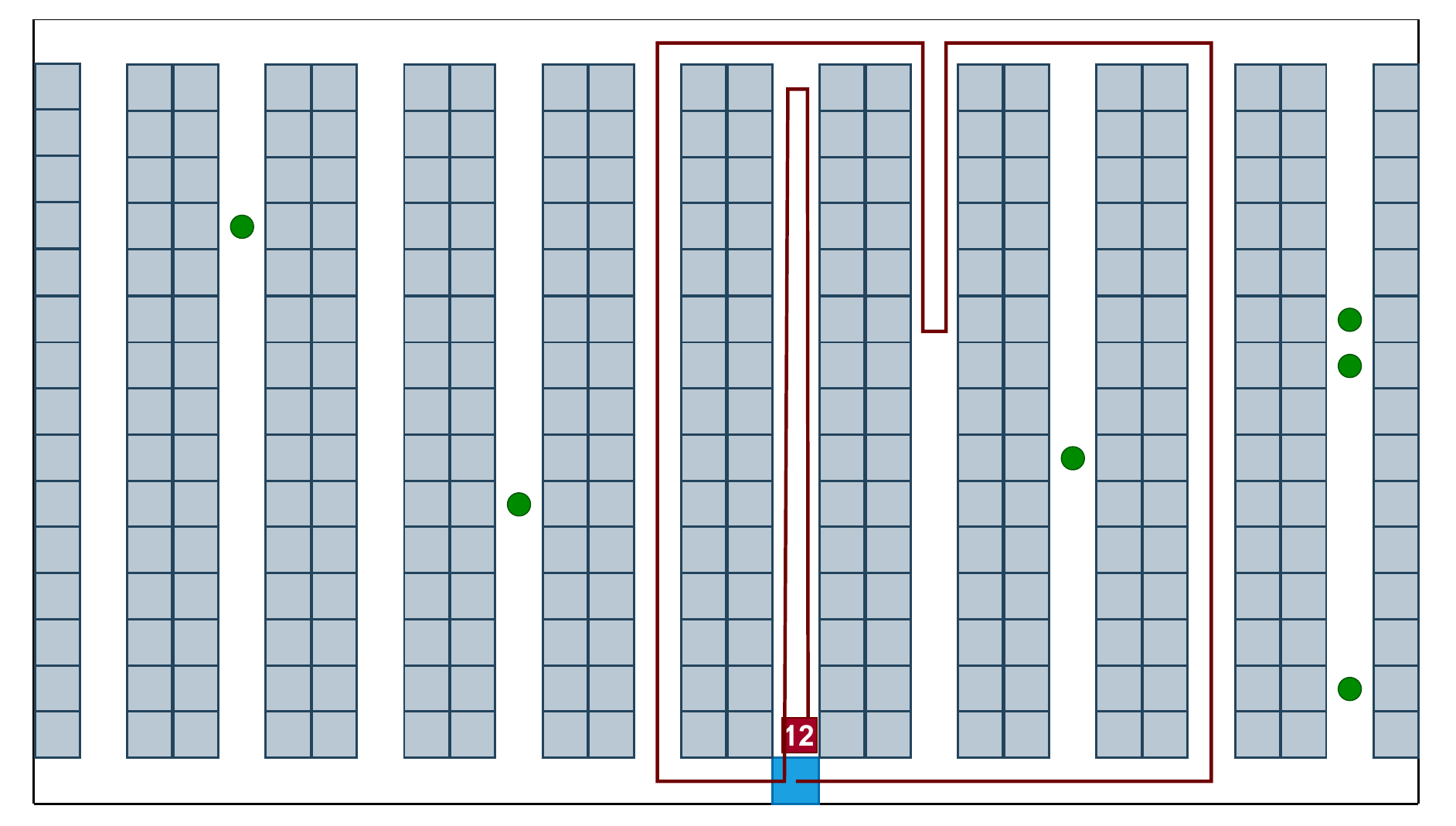}
        \caption{Time (sec.): 766}
        \label{subfig:766}
    \end{subfigure}
    \hfill
    \begin{subfigure}[b]{0.49\textwidth}
        \centering
        \includegraphics[width=\textwidth]{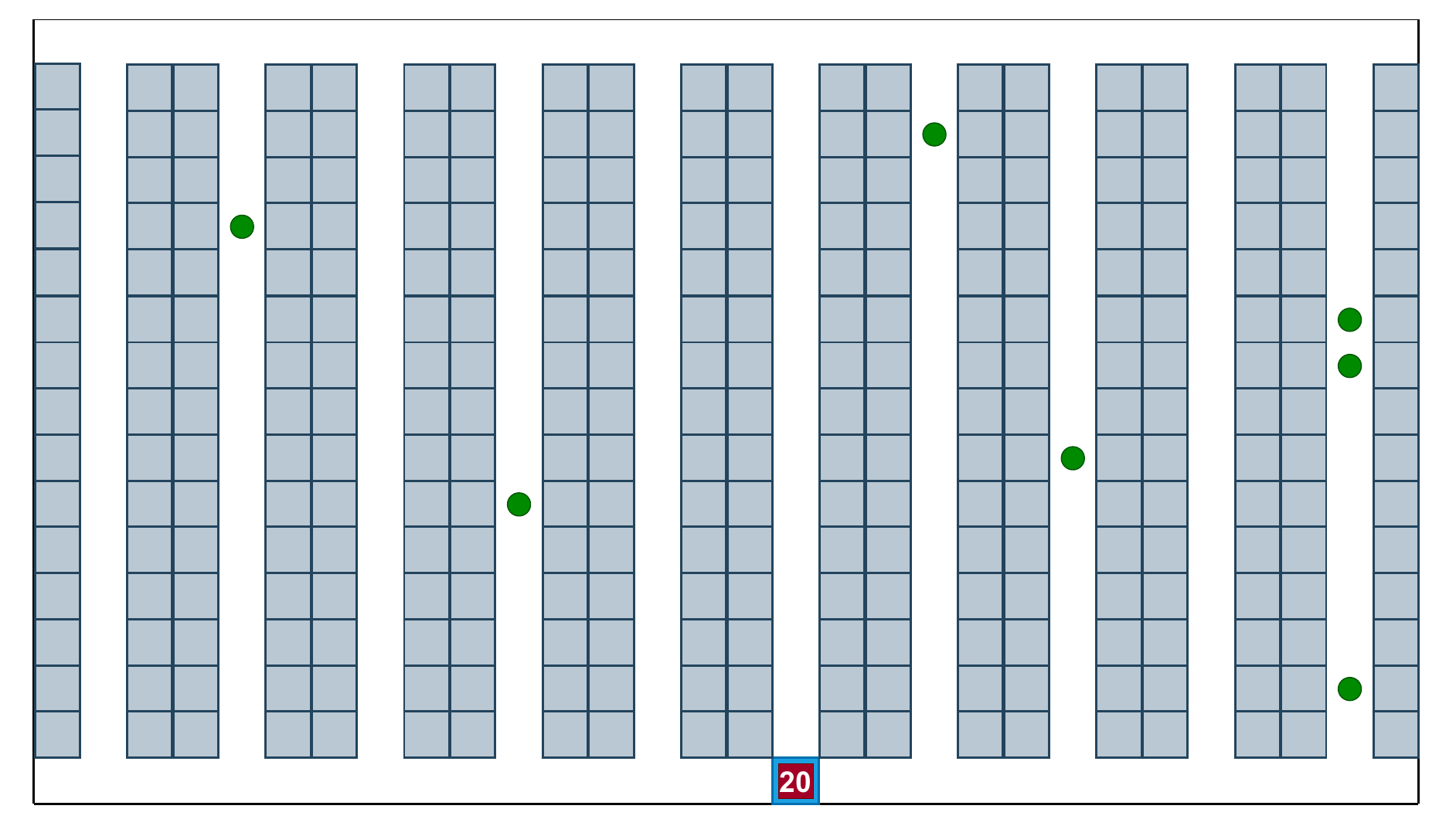}
        \caption{Time (sec.): 775}
        \label{subfig:775}
    \end{subfigure}
    \caption{Example demonstration of a pick cycle executed by the learned DRL agent with  $\lambda=0.06$ and $\alpha = 1$}
    \label{fig:illustration}
\end{figure}

\paragraph{Illustration of Real-time Order Picking} Interested readers may refer to a YouTube animation video (\href{https://youtu.be/LumE7YE-JMo}{https://youtu.be/LumE7YE-JMo}) that illustrates the real-time progress of the DRL agent. However, in this section, we present a few representative snapshots to highlight key aspects of the agent’s behavior over time. Specifically, figure~\ref{fig:illustration} illustrates a pick cycle recommended by our trained DRL agent for a specific scenario with $\lambda=0.06$ and $\alpha=1$. In this context, the term `time' refers to the elapsed time since the start of the work shift (not the computational time). Accordingly, the picker begins this cycle 661 seconds into the work shift with a full capacity of $K=20$. The blue square represents the depot, while the red square denotes the picker with its remaining capacity highlighted. Since the horizontal distance within aisles is negligible, orders within storage locations on both sides are shown as green dots inside the aisles. At 673 seconds, as depicted in Figure~\ref{subfig:673}, the order picker arrives at aisle 9 and retrieves an item from the second row. Simultaneously, a new demand arises for an item on the third row in the same aisle. The picker efficiently fulfills all three pick requests in aisle 9 before moving sequentially to aisles 7 and 5, as indicated by the red line. Figure~\ref{subfig:735} illustrates the picker carrying five items with a remaining capacity of 15 units as it proceeds towards aisle 6 to pick up the next set of requested items. At this point, two new orders appear at position 13 in aisle 6. The DRL agent’s demand forecasting guided the picker to this aisle. Figure~\ref{subfig:766} shows the picker completing pick-ups in aisle 6 before returning to the depot via the same path. In Figure~\ref{subfig:775}, the picker concludes the drop-offs at the depot, where each drop-off takes one second. This marks the end of the pick cycle, and the picker is prepared for the next assignment. Note that new orders have arrived at the warehouse at various time steps during the picker’s traversal. Throughout the pick cycle, the DRL agent dynamically prioritizes orders to minimize the average order fulfillment time.

\subsection{Numerical Results}

The results of our experiments are summarized in Tables \ref{tab:ATDO}, \ref{tab:AOCT}, and \ref{tab:PUO}. 
Table~\ref{tab:ATDO} presents the average total distance traveled (ATDO) resulting from each model decisions for each instance type.
The corresponding results for average order completion time (AOCT) and percentage of unfulfilled orders (PUO) are shown in Tables \ref{tab:AOCT} and \ref{tab:PUO}, respectively. Across these tables, the shaded cells highlight either the most desirable result or the most effective DRL model. 
Specifically, the trained model demonstrating the best average order completion time for each order arrival rate is highlighted. Next, we analyse the results and make a few observations. 

The first observation is that, unsurprisingly, Baseline 1 results in the shortest distance traveled due to its strategy of waiting for the pick list to reach pick device capacity. This was expected, as by design, this approach allows the picker to follow the shortest possible path and move uninterrupted even when new customer orders arrive. However, as previously mentioned, the most crucial metric for optimizing warehouse operations is not the shortest distance traveled but the order completion times. Unfortunately, Baseline 1 has the worst performance in terms of average order completion times across all arrival rates and, consequently, the highest percentage of unfilled orders. Therefore, among the baselines, the variations of IRA—namely Baselines 2 through 5—are certainly better choices than Baseline 1. 
Notably, heuristic policies such as S-shape and Largest-gap consistently perform poorly compared to their optimal distance-routing counterpart, Baseline 5, across all performance metrics. This outcome is unsurprising, as the routing decisions of these heuristics are suboptimal compared to their corresponding exact baselines, which rely on solving the TSP to optimality at all times.


\begin{table}[]
\scriptsize
    \centering
    \caption{Performance Comparison Summary: Average Travel Distance per Order (ATDO)}
    \label{tab:ATDO}
    \begin{tabular}{cc rrrrrrrrr}
    \toprule
\multicolumn{2}{c}{DRL Model}                & \multicolumn{9}{c}{ATDO (meters) w.r.t. $\lambda$} \\
\cmidrule(lr){1-2} \cmidrule(lr){3-11}
$\lambda$ & $\alpha$ & \multicolumn{1}{c}{0.01}   & \multicolumn{1}{c}{0.02}   & \multicolumn{1}{c}{0.03}   & \multicolumn{1}{c}{0.04}   & \multicolumn{1}{c}{0.05}   & \multicolumn{1}{c}{0.06}   & \multicolumn{1}{c}{0.07}  & \multicolumn{1}{c}{0.08}  & \multicolumn{1}{c}{0.09}  \\
 \midrule
 & 0.1         & \cellcolor[gray]{0.9}37.73 & \cellcolor[gray]{0.9}55.83 & \cellcolor[gray]{0.9}86.42 & \cellcolor[gray]{0.9}111.15 & 109.81  & 95.24   & 80.81   & 65.61   & 48.50   \\
 & 0.5         & \cellcolor[gray]{0.9}37.15 & \cellcolor[gray]{0.9}45.28 & \cellcolor[gray]{0.9}57.43 & \cellcolor[gray]{0.9}63.91  & 66.99   & 69.61   & 68.18   & 61.40   & 51.53   \\
\multirow{-3}{*}{0.02} & 1.0         & \cellcolor[gray]{0.9}33.85 & \cellcolor[gray]{0.9}40.31 & \cellcolor[gray]{0.9}48.70 & \cellcolor[gray]{0.9}59.89  & 73.69   & 86.12   & 81.71   & 70.74   & 56.97   \\
\hline
 & 0.1         & 61.93   & 79.65   & 87.51   & 88.83    & \cellcolor[gray]{0.9}84.83 & \cellcolor[gray]{0.9}77.18 & 69.63   & 63.03   & 55.31   \\
 & 0.5         & 81.60   & 86.40   & 81.57   & 78.43    & \cellcolor[gray]{0.9}74.42 & \cellcolor[gray]{0.9}70.68 & 67.46   & 65.37   & 57.75   \\
\multirow{-3}{*}{0.04} & 1.0         & 72.17   & 75.91   & 71.71   & 66.80    & \cellcolor[gray]{0.9}63.30 & \cellcolor[gray]{0.9}60.74 & 59.54   & 57.03   & 48.58   \\
\hline
 & 0.1         & 65.78   & 106.23  & 106.92  & 88.80    & \cellcolor[gray]{0.9}75.32 & \cellcolor[gray]{0.9}61.26 & \cellcolor[gray]{0.9}53.19 & 46.19   & 52.29   \\
 & 0.5         & 61.51   & 93.21   & 95.26   & 81.43    & \cellcolor[gray]{0.9}69.21 & \cellcolor[gray]{0.9}60.15 & \cellcolor[gray]{0.9}56.10 & 49.52   & 46.36   \\
\multirow{-3}{*}{0.06} & 1.0         & 117.30  & 107.46  & 87.08   & 75.31    & \cellcolor[gray]{0.9}65.38 & \cellcolor[gray]{0.9}59.03 & \cellcolor[gray]{0.9}55.66 & 51.27   & 49.72   \\
\hline
 & 0.1         & 71.71   & 145.70  & 198.54  & 177.00   & 143.14  & 103.99  & 75.88   & \cellcolor[gray]{0.9}57.51 & \cellcolor[gray]{0.9}47.08 \\
 & 0.5         & 126.68  & 132.04  & 109.68  & 83.98    & 69.97   & 59.70   & 54.62   & \cellcolor[gray]{0.9}49.27 & \cellcolor[gray]{0.9}42.26 \\
\multirow{-3}{*}{0.08} & 1.0         & 70.31   & 126.53  & 141.88  & 112.20   & 87.12   & 67.17   & 54.66   & \cellcolor[gray]{0.9}46.57 & \cellcolor[gray]{0.9}40.00 \\
\hline
\hline
\multicolumn{2}{c}{Baseline 1}  & \cellcolor[gray]{0.9}8.18   & \cellcolor[gray]{0.9}8.14   & \cellcolor[gray]{0.9}8.20   & \cellcolor[gray]{0.9}8.18   & \cellcolor[gray]{0.9}8.14   & \cellcolor[gray]{0.9}8.09   & \cellcolor[gray]{0.9}8.19  & \cellcolor[gray]{0.9}8.17  & \cellcolor[gray]{0.9}8.19  \\
\multicolumn{2}{c}{Baseline 2}   & 16.69  & 16.43  & 15.60  & 14.72  & 13.05  & 10.90  & 8.84  & 8.28  & 8.23  \\
\multicolumn{2}{c}{Baseline 3}  & 16.78  & 16.59  & 15.71  & 14.68  & 13.01  & 10.88  & 8.84  & 8.28  & 8.23  \\
\multicolumn{2}{c}{Baseline 4}   & 27.98  & 24.65  & 21.08  & 17.56  & 14.39  & 11.18  & 8.88  & 8.28  & 8.26  \\
\multicolumn{2}{c}{Baseline 5}  & 27.74  & 24.37  & 20.86  & 17.45  & 14.35  & 11.15  & 8.88  & 8.28  & 8.26 \\
\multicolumn{2}{c}{\textcolor{black}{S-shape}} & \textcolor{black}{29.86} & \textcolor{black}{28.20} & \textcolor{black}{24.55} & \textcolor{black}{19.21} & \textcolor{black}{14.96} & \textcolor{black}{12.63} & \textcolor{black}{12.58} & \textcolor{black}{12.60} & \textcolor{black}{12.55} \\
\multicolumn{2}{c}{\textcolor{black}{Largest-gap}} & \textcolor{black}{29.95} & \textcolor{black}{27.70} & \textcolor{black}{23.96} & \textcolor{black}{19.00} & \textcolor{black}{14.86} & \textcolor{black}{11.28} & \textcolor{black}{10.33} & \textcolor{black}{10.26} & \textcolor{black}{10.24} \\

\bottomrule 
\end{tabular}
\end{table}

Further analysis of baseline models also reveals that allowing pickers to re-route within cross-aisles enhances order completion times when order arrivals are frequent, i.e.,  $\lambda \geq 0.05$.
However, the intervention routing strategy should be avoided within cross-aisles when order arrivals are less frequent ($\lambda < 0.05$).
This is evident when comparing Baselines 2 and 3, or Baselines 4 and 5, in Table~\ref{tab:AOCT}. 
Notably, under scenarios with higher order arrival rates, all baseline models lead to a substantial number of unfulfilled orders within the work shift. 
For example, when $\lambda=0.9$, approximately 18\% of orders remain unfulfilled. 
An exception occurs when $\lambda \leq 0.02$, where Baseline 4 surpasses our trained DRL models. 
In such cases of low arrival rates, a simple policy of shortest routing with intervention can be effective.

\begin{table}[]
\scriptsize
    \centering
    \caption{Performance Comparison Summary: Average Order Completion Time (AOCT)}
    \label{tab:AOCT}
    \begin{tabular}{cc rrrrrrrrr}
    \toprule
\multicolumn{2}{c}{DRL Model}                & \multicolumn{9}{c}{AOCT (seconds) w.r.t. $\lambda$} \\
\cmidrule(lr){1-2} \cmidrule(lr){3-11}
$\lambda$ & $\alpha$ & \multicolumn{1}{c}{0.01}   & \multicolumn{1}{c}{0.02}   & \multicolumn{1}{c}{0.03}   & \multicolumn{1}{c}{0.04}   & \multicolumn{1}{c}{0.05}   & \multicolumn{1}{c}{0.06}   & \multicolumn{1}{c}{0.07}  & \multicolumn{1}{c}{0.08}  & \multicolumn{1}{c}{0.09}  \\
 \midrule
 & 0.1                  & \cellcolor[gray]{0.9}89.9 & \cellcolor[gray]{0.9}122.0 & \cellcolor[gray]{0.9}182.9 & \cellcolor[gray]{0.9}268.4 & 307.6   & 368.4   & 461.3   & 550.2   & 764.5   \\
 & 0.5                  & \cellcolor[gray]{0.9}81.5 & \cellcolor[gray]{0.9}95.1  & \cellcolor[gray]{0.9}129.6 & \cellcolor[gray]{0.9}170.8 & 219.6   & 300.9   & 392.6   & 479.7   & 644.5   \\
\multirow{-3}{*}{0.02} & 1.0                  & \cellcolor[gray]{0.9}67.4 & \cellcolor[gray]{0.9}87.2  & \cellcolor[gray]{0.9}118.1 & \cellcolor[gray]{0.9}171.6 & 239.8   & 363.4   & 512.2   & 676.3   & 844.9   \\
\hline
 & 0.1                  & 150.5  & 160.9   & 178.8   & 208.3   & \cellcolor[gray]{0.9}238.2 & \cellcolor[gray]{0.9}278.2 & 340.8   & 460.5   & 740.1   \\
 & 0.5                  & 226.5  & 181.2   & 174.8   & 196.5   & \cellcolor[gray]{0.9}222.0 & \cellcolor[gray]{0.9}266.4 & 336.0   & 456.5   & 716.8   \\
\multirow{-3}{*}{0.04} & 1.0                  & 180.1  & 156.8   & 162.4   & 185.5   & \cellcolor[gray]{0.9}210.2 & \cellcolor[gray]{0.9}255.9 & 330.5   & 436.6   & 620.3   \\
\hline
 & 0.1                  & 136.0  & 192.6   & 218.8   & 227.3   & \cellcolor[gray]{0.9}236.3 & \cellcolor[gray]{0.9}262.3 & \cellcolor[gray]{0.9}305.1 & 402.0   & 809.6   \\
 & 0.5                  & 146.8  & 186.3   & 208.8   & 218.3   & \cellcolor[gray]{0.9}230.3 & \cellcolor[gray]{0.9}256.1 & \cellcolor[gray]{0.9}304.2 & 382.1   & 676.8   \\
\multirow{-3}{*}{0.06} & 1.0                  & 277.6  & 218.1   & 202.8   & 216.7   & \cellcolor[gray]{0.9}227.1 & \cellcolor[gray]{0.9}260.9 & \cellcolor[gray]{0.9}323.1 & 403.1   & 640.2   \\
\hline
 & 0.1                  & 175.1  & 332.8   & 460.3   & 427.0   & 400.5   & 378.9   & 378.3   & \cellcolor[gray]{0.9}403.7 & \cellcolor[gray]{0.9}513.1 \\
 & 0.5                  & 309.0  & 292.0   & 294.7   & 284.7   & 288.7   & 306.0   & 336.0   & \cellcolor[gray]{0.9}383.5 & \cellcolor[gray]{0.9}536.7 \\
\multirow{-3}{*}{0.08} & 1.0                  & 158.2  & 252.5   & 313.0   & 304.4   & 291.8   & 295.0   & 318.5   & \cellcolor[gray]{0.9}369.0 & \cellcolor[gray]{0.9}556.3 \\ 
\hline
\hline
\multicolumn{2}{c}{Baseline 1}  & 1,217.1 & 751.6 & 593.0 & 512.1 & 471.7 & 445.1 & 560.0 & 1,526.2 & 2,809.3 \\
\multicolumn{2}{c}{Baseline 2}   & 292.6   & 226.3 & 238.0 & 270.1 & 289.2 & 313.1 & 461.5 & 1,429.5 & 2,715.8 \\
\multicolumn{2}{c}{Baseline 3}  & 294.5   & 229.0 & 242.9 & 270.2 & 286.8 & 311.6 & 460.2 & 1,425.6 & 2,710.1 \\
\multicolumn{2}{c}{Baseline 4}   & \cellcolor[gray]{0.9}52.2    & \cellcolor[gray]{0.9}82.1  & 141.4 & 222.0 & 271.2 & 308.7 & 462.9 & 1,419.8 & 2,712.7 \\
\multicolumn{2}{c}{Baseline 5}  & 54.6    & 87.1  & 148.3 & 225.2 & 271.0 & 308.8 & 461.2 & 1,418.4 & 2,711.8 \\
\multicolumn{2}{c}{\textcolor{black}{S-shape}} & \textcolor{black}{65.0} & \textcolor{black}{154.4} & \textcolor{black}{304.2} & \textcolor{black}{409.8} & \textcolor{black}{484.1} & \textcolor{black}{1,546.0} & \textcolor{black}{3,283.9} & \textcolor{black}{4,492.8} & \textcolor{black}{5,432.5} \\
\multicolumn{2}{c}{\textcolor{black}{Largest-gap}} & \textcolor{black}{65.0} & \textcolor{black}{138.1} & \textcolor{black}{274.3} & \textcolor{black}{373.1} & \textcolor{black}{401.6} & \textcolor{black}{494.3} & \textcolor{black}{1,707.1} & \textcolor{black}{3,029.1} & \textcolor{black}{4,154.8} \\
\bottomrule 
\end{tabular}
\end{table}

\begin{table}[]
    \centering
    \scriptsize
    \caption{Performance Comparison Summary: Percentage of Unfulfilled Orders (PUO)}
    \label{tab:PUO}
    \begin{tabular}{cc rrrrrrrrr}
    \toprule
\multicolumn{2}{c}{DRL Model}                & \multicolumn{9}{c}{PUO (\%) w.r.t. $\lambda$} \\
\cmidrule(lr){1-2} \cmidrule(lr){3-11}
$\lambda$ & $\alpha$ & \multicolumn{1}{c}{0.01}   & \multicolumn{1}{c}{0.02}   & \multicolumn{1}{c}{0.03}   & \multicolumn{1}{c}{0.04}   & \multicolumn{1}{c}{0.05}   & \multicolumn{1}{c}{0.06}   & \multicolumn{1}{c}{0.07}  & \multicolumn{1}{c}{0.08}  & \multicolumn{1}{c}{0.09}  \\
 \midrule
 & 0.1 & \cellcolor[gray]{0.9}0.46 & \cellcolor[gray]{0.9}0.38 & \cellcolor[gray]{0.9}0.36 & \cellcolor[gray]{0.9}1.22 & 1.11   & 1.24   & 1.26   & 2.13   & 2.80   \\
 & 0.5 & \cellcolor[gray]{0.9}0.40 & \cellcolor[gray]{0.9}0.27 & \cellcolor[gray]{0.9}0.45 & \cellcolor[gray]{0.9}0.54 & 0.66   & 0.91   & 1.21   & 1.82   & 2.16   \\
\multirow{-3}{*}{0.02} & 1.0 & \cellcolor[gray]{0.9}0.64 & \cellcolor[gray]{0.9}0.36 & \cellcolor[gray]{0.9}0.30 & \cellcolor[gray]{0.9}0.60 & 0.85   & 1.20   & 1.72   & 2.72   & 2.99   \\
\hline
 & 0.1 & 0.63   & 0.54   & 0.60   & 0.93   & \cellcolor[gray]{0.9}0.93 & \cellcolor[gray]{0.9}0.71 & 1.01   & 1.62   & 3.02   \\
 & 0.5 & 1.23   & 0.67   & 0.50   & 0.75   & \cellcolor[gray]{0.9}0.69 & \cellcolor[gray]{0.9}0.72 & 1.00   & 1.77   & 2.78   \\
\multirow{-3}{*}{0.04} & 1.0 & 0.63   & 0.67   & 0.60   & 0.74   & \cellcolor[gray]{0.9}0.56 & \cellcolor[gray]{0.9}0.79 & 0.92   & 1.69   & 2.19   \\
\hline
 & 0.1 & 0.76   & 0.57   & 0.69   & 0.72   & \cellcolor[gray]{0.9}0.75 & \cellcolor[gray]{0.9}0.80 & \cellcolor[gray]{0.9}1.01 & 1.47   & 5.62   \\
 & 0.5 & 0.84   & 0.65   & 0.62   & 0.93   & \cellcolor[gray]{0.9}0.74 & \cellcolor[gray]{0.9}0.71 & \cellcolor[gray]{0.9}0.95 & 1.40   & 3.47   \\
\multirow{-3}{*}{0.06} & 1.0 & 0.80   & 0.93   & 0.57   & 0.73   & \cellcolor[gray]{0.9}0.69 & \cellcolor[gray]{0.9}0.69 & \cellcolor[gray]{0.9}1.01 & 1.54   & 2.66   \\
\hline
 & 0.1 & 0.81   & 0.89   & 1.21   & 1.82   & 1.31   & 1.27   & 1.36   & \cellcolor[gray]{0.9}1.36 & \cellcolor[gray]{0.9}1.78 \\
 & 0.5 & 1.87   & 0.86   & 1.12   & 1.11   & 0.86   & 0.95   & 1.17   & \cellcolor[gray]{0.9}1.27 & \cellcolor[gray]{0.9}1.89 \\
\multirow{-3}{*}{0.08} & 1.0 & 0.66   & 0.77   & 1.08   & 1.25   & 0.96   & 1.03   & 0.97   & \cellcolor[gray]{0.9}1.25 & \cellcolor[gray]{0.9}2.27
 \\
\hline
\hline
\multicolumn{2}{c}{Baseline 1}   & 5.02 & 1.89 & 2.26 & 1.48 & 1.43 & 1.22 & 1.65 & 9.24 & 18.55 \\
\multicolumn{2}{c}{Baseline 2}   & 0.92 & 0.57 & 0.58 & 1.02 & 0.71 & 0.89 & 1.15 & 9.06 & 18.30 \\
\multicolumn{2}{c}{Baseline 3}  & 0.81 & 0.55 & 0.68 & 1.06 & 0.83 & 0.88 & 1.14 & 9.06 & 18.30 \\
\multicolumn{2}{c}{Baseline 4}   & \cellcolor[gray]{0.9}0.07 & \cellcolor[gray]{0.9}0.25 & 0.32 & 0.94 & 0.69 & 1.17 & 1.18 & 8.86 & 18.30 \\
\multicolumn{2}{c}{Baseline 5}  & 0.10 & \cellcolor[gray]{0.9}0.25 & 0.34 & 0.99 & 0.70 & 0.93 & 1.13 & 8.86 & 18.30 \\
\multicolumn{2}{c}{\textcolor{black}{S-shape}} & \textcolor{black}{0.32} & \textcolor{black}{0.54} & \textcolor{black}{0.77} & \textcolor{black}{1.76} & \textcolor{black}{1.50} & \textcolor{black}{8.64} & \textcolor{black}{20.82} & \textcolor{black}{30.06} & \textcolor{black}{37.19} \\
\multicolumn{2}{c}{\textcolor{black}{Largest-gap}} & \textcolor{black}{0.32} & \textcolor{black}{0.40} & \textcolor{black}{0.62} & \textcolor{black}{1.71} & \textcolor{black}{1.21} & \textcolor{black}{1.49} & \textcolor{black}{9.99} & \textcolor{black}{20.00} & \textcolor{black}{28.27} \\

\bottomrule 
\end{tabular}
\end{table}

Comparing the baselines with our trained DRL models demonstrates that our DRL model significantly enhances both order throughput and fulfillment rates, even in scenarios not encountered during training. Based on the results from the trained DRL models, we recommend the following policies for training DRL models: Specifically, use training data that follows Poisson arrival rate parameters aligned with actual order arrival rates, as follows:

\begin{itemize}
    \item Train with $\lambda = 0.02$: For observed order arrival rates between 0.01 and 0.04
    \item Train with $\lambda = 0.06$: For observed order arrival rates between 0.05 and 0.07
    \item Train with $\lambda = 0.08$: For observed order arrival rates between 0.08 and 0.09
\end{itemize}

While the learned policies may lead to increased travel distances within the warehouses, they result in significant time savings overall. 
For example, in scenarios with $\lambda=0.09$, the recommended DRL model's decision rules reduce AOCT by about 428\% compared to baseline models. 
Moreover, our model drastically reduces the percentage of unfulfilled orders (PUO) to approximately 2\%, compared to 18.3\% for baselines during the 8-hour work shifts. 
Our results demonstrate that demand forecasting and order picking using reinforcement learning are highly effective strategies, particularly in scenarios with frequent order arrivals, i.e., $\lambda > 0.02$.

Finally, our experiments demonstrate the remarkable robustness and transferability of our trained models across diverse testing scenarios with varying arrival rates.
For each $\lambda$ used in training, we explored three values of $\alpha$, a parameter indirectly controlling the trade-off between prioritizing travel distance ($\alpha =0$) and order completion times ($\alpha = 1$).
While `$\alpha = 1$' performs well when training and testing data are closely aligned, we find that `$\alpha = 0.5$' provides superior robustness when real arrival rates may fluctuate slightly.
Additionally, in test instances with higher arrival rates, all our trained models significantly outperformed benchmark algorithms in terms of AOCT and PUO.
For example, in test instances with an arrival rate of $\lambda = 0.08$, the worst-performing trained model still achieved an AOCT of 676.3 seconds, a significant improvement over the best baseline's AOCT of 1418.4 seconds. 
These results underscore the real-world applicability of our proposed framework, as it maintains exceptional performance even under varying order arrival rates, making it well-suited for real-time deployment.

\begin{figure}[ht]
    \centering
    \begin{subfigure}[b]{0.49\textwidth}
        \centering
        \includegraphics[width=\textwidth]{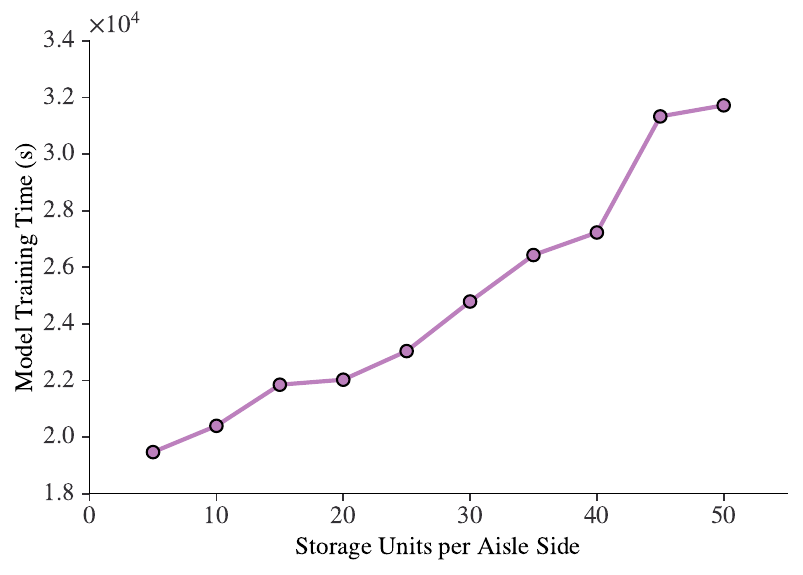}
        \caption{Effect of the number of storage locations per aisle when the number of aisles is 10}
        \label{fig:runtime_rows}
    \end{subfigure}
    \hfill
    \begin{subfigure}[b]{0.49\textwidth}
        \centering
        \includegraphics[width=\textwidth]{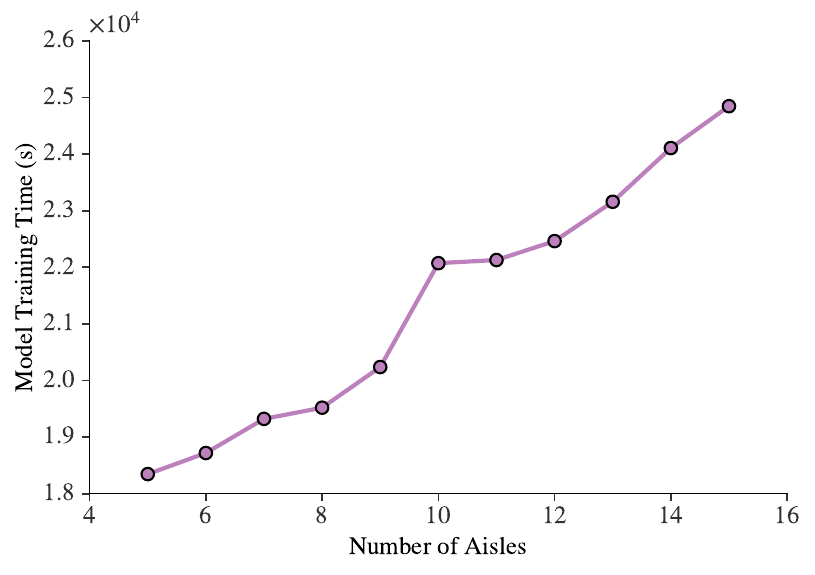} 
        \caption{Effect of the number of aisles when the number of storage locations per aisle is 15}
        \label{fig:runtime_aisles}
    \end{subfigure}
    \caption{Training time of the DRL model for $\lambda = 0.08$ and $\alpha = 1$ across different warehouse configurations}
    \label{fig:runtimes}
\end{figure}

\section{Final Remarks} \label{sec:conclusions}

In this research, we proposed a DRL model designed to solve the dynamic order picking problem within a single-block warehouse layout served by an autonomous picking device. 
Our results demonstrate that the model's ability to forecast customer demands and adaptively learn efficient order picking policies leads to significant enhancements in overall warehouse operations. 
We benchmarked the proposed model against optimal static algorithms and widely used heuristics designed for picking distance optimization, and showed that our approach achieves a substantial reduction in average order completion times while fulfilling a significantly higher number of orders during the work shift. 
Furthermore, additional experiments validated the robustness of the model, as evidenced by its strong performance in out-of-sample test instances with varying order arrival rates. 
Our findings highlighted the adaptability and reliability of the DRL framework and established it as a promising method for advancing warehouse automation.

As mentioned, the proposed method has been shown to be suitable for single-block warehouses with a single autonomous picking device. However, even within its domain of applicability, it is important to highlight that some computational challenges may arise depending on the structure of the warehouse. Specifically, we demonstrated the effectiveness of our approach on warehouses with 10 aisles and 15 storage locations per side of each aisle, which is a common assumption in the existing body of literature. If the size of the warehouse changes, the proposed method remains effective but is expected to require significantly more computational power during the training phase. This is illustrated in Figure~\ref{fig:runtimes}, where the computational time required to train the DRL model under varying warehouse configurations is presented with $\lambda = 0.08$ and $\alpha = 1$. Figure~\ref{fig:runtime_rows} examines the impact of increasing the number of storage locations per aisle while keeping the number of aisles fixed at 10. Conversely, Figure~\ref{fig:runtime_aisles} explores the effect of increasing the number of aisles, with the number of storage locations per aisle side fixed at 15. In both scenarios, the results reveal a clear trend: training time increases substantially as warehouse size grows. 
Despite this computational cost, the DRL model demonstrates significant potential for practical problem-solving. Most warehouse configurations remain fixed once determined, making the one-time training overhead manageable. 
To further alleviate this computational burden, parallel computing techniques could be leveraged to distribute training across multiple GPUs or processing units.

It is also worth mentioning that this work provides a foundation for future research to explore more complex warehouse environments. 
While this study focuses on a single-block warehouse configuration, extending the model to multi-block layouts offers the opportunity to address the greater complexities of real-world operations. Similarly, extending the current scenario from a single autonomous picking device to one involving multiple coordinated devices, whether fully autonomous or incorporating human-robot interactions, presents a compelling avenue for further exploration. 
In both of these research avenues, employing parallelization and exploring suitable neural network architectures and reward functions will emerge as the primary challenges.

\bibliographystyle{ormsv080} 
\bibliography{ref}

\newpage
\appendix

\section{Rationale for Feature-Based State Representation}
\label{app:rationale}

A significant challenge in dynamic order-picking tasks is managing the extremely high-dimensional raw data, such as exact picker coordinates and detailed item-level positions for every storage location. To overcome this, we propose an aggregated feature-based state representation that retains all decision-relevant information while avoiding unnecessary fine-grained spatial details. Our representation captures essential routing information, picker constraints, and order demands through targeted aggregation.

We aggregate orders at the aisle level by separately considering upward and downward traversal directions within each aisle. Specifically, each aisle is represented by two aggregated indices corresponding to the potential rewards of traversing upward and downward. Potential rewards are calculated dynamically, considering both the demand (number of items) and their distance from the picker's current position. Items closer to the picker have higher reward contributions, thereby guiding the agent to prioritize more immediately accessible orders.

The picker's state is represented concisely through three main components:
\begin{itemize}
    \item \textbf{Horizontal Position} $(S^H_t)$: Encodes the picker's horizontal location as either in the back cross-aisle $(S^H_t = -1)$, within a pick aisle $(S^H_t = 0)$, or in the front cross-aisle $(S^H_t = 1)$. This abstract encoding preserves essential movement constraints without the complexity of precise coordinates.

    \item \textbf{Vertical Position} $(S^{V_1}_t, S^{V_2}_t)$: Each aisle $n$ is associated with two distinct indices, $2n-1$ and $2n$, which align with upward and downward traversal rewards, respectively. When the picker is positioned in aisle $n$, we set $(S^{V_1}_t, S^{V_2}_t) = (2n-1, 2n)$. This indexing scheme directly aligns picker positions with the order state representation, enabling efficient computation of aisle-specific potential rewards.

    \item \textbf{Remaining Capacity} $(S^C_t)$: Indicates the available capacity of the picker, ranging from fully loaded $(0)$ to completely empty $(K)$. Including capacity explicitly allows the agent to plan efficient routes, avoid unnecessary visits to aisles when capacity is insufficient, and effectively schedule depot returns. Representing capacity with a single scalar prevents unnecessary dimensionality growth while maintaining critical constraints.
\end{itemize}

The aggregate order state is encoded as:
\[
\boldsymbol{S^o_t} = (S^1_t, S^2_t, \ldots, S^{2N-1}_t, S^{2N}_t),
\]
where each aisle $n$ corresponds to two components, $(S^{2n-1}_t, S^{2n}_t)$. These values are dynamically computed, capturing the combined item demand weighted by distance, separately for upward and downward traversals within each aisle. As picker location or order demand changes, these reward potentials are continuously updated, ensuring the representation reflects the current system state without explicit spatial detail.

By adopting this feature-based representation, we effectively eliminate redundant fine-grained spatial data, significantly reducing the dimensionality of the state space. The proposed representation retains critical information required for efficient routing decisions---such as aisle-specific potential rewards and capacity constraints---without overwhelming the agent with irrelevant coordinate-level detail. Consequently, this facilitates faster policy learning, improved interpretability, and greater robustness against overfitting to unnecessary spatial complexity.

\section{Analysis of the DRL Model with Different Problem Settings} 
\label{app:analysis}

Our DRL model assumes that the picker can initiate movement for order picking from the depot as long as at least one order is pending.
However, studies such as \citet{ratliff1983order} and \citet{lu2016algorithm} establish a higher bound on the initial pick size for picker routing.
Furthermore, these studies assume narrow cross-aisles, preventing rerouting through them. This limitation ensures that the picker does not move back in the opposite direction without entering an aisle.
In contrast, our approach permits cross-aisle rerouting in the primary experiments, allowing for greater flexibility in picker movement. In this section, we incorporate a combination of these restrictive assumptions into the DRL agent’s environment and evaluate its performance after training under these conditions.

Table \ref{tab:Models} summarizes the various DRL models used for comparison. In this section, the parameter $\alpha$ is set to 1. Additionally, we emphasize that the main manuscript includes an analysis of the arrival rate on which the DRL agent should be trained to perform well across a range of arrival rates during testing. However, in this section, we do not repeat that analysis. Instead, we directly apply the insights gained from the main experiments in the manuscript. Specifically, the arrival rate ($\lambda$) used for training each model is chosen based on the best observed performance in terms of average order completion times, as reported in the manuscript’s primary experiments (i.e., the gray-highlighted sections in Table~\ref{tab:AOCT}).
Specifically, for $\lambda\in\{0.01, 0.02, 0.03, 0.04\}$, the DRL model trained on an arrival rate of 0.02 is utilized. For $\lambda = 0.05$, the model trained on an arrival rate of 0.04 is selected. Similarly, for $\lambda\in\{0.06, 0.07\}$, the model trained on an arrival rate of 0.06 is applied, and for $\lambda\in\{0.08, 0.09\}$, the model trained on an arrival rate of 0.08 is used. 
This approach ensures the selection of models best suited to optimize performance under varying arrival rate conditions, allowing for a focused analysis of the maximum achievable performance across different variations of the DRL model highlighted in Table~\ref{tab:Models}.
\begin{table}[!htbp]
    \centering
    \scriptsize
    \caption{Variations of the DRL Model Used for Sensitivity Analysis}
    \label{tab:Models}
    \begin{tabular}{ccc}
    \toprule
    Model      & Initial Pick Size  & Cross-aisle Re-routing \\
    \midrule
    DRL 1 & 20 & No                   \\
    DRL 2 & 20 & Yes                   \\
    DRL 3 & 5  & No                     \\
    DRL 4 & 5  & Yes                    \\
    DRL 5 & 1  & No                     \\
    DRL 6 & 1  & Yes                    \\
    \bottomrule
    \end{tabular}
\end{table}

\begin{table}[]
    \centering
    \scriptsize
    \caption{Sensitivity Analysis: Average Travel Distance per Order (ATDO)}
    \label{tab:app_ATDO}
    \begin{tabular}{cc rrrrrrrrr}
    \toprule
\multicolumn{2}{c}{}                & \multicolumn{9}{c}{ATDO (meters) w.r.t. $\lambda$} \\
 \cmidrule(lr){3-11}
\multicolumn{2}{c}{Model} & \multicolumn{1}{c}{0.01}   & \multicolumn{1}{c}{0.02}   & \multicolumn{1}{c}{0.03}   & \multicolumn{1}{c}{0.04}   & \multicolumn{1}{c}{0.05}   & \multicolumn{1}{c}{0.06}   & \multicolumn{1}{c}{0.07}  & \multicolumn{1}{c}{0.08}  & \multicolumn{1}{c}{0.09}  \\
 \midrule
 \multicolumn{2}{c}{DRL 1} & 49.85 & 52.15 & 56.58 & 67.46 & 54.85 & 60.89 & 48.76 & 69.61 & 53.54 \\  \multicolumn{2}{c}{DRL 2} & 51.46 & 54.72 & 60.03 & 69.91 & 58.3 & 55.82 & 51.89 & 45.79 & 38.88 \\  \multicolumn{2}{c}{DRL 3} & 48.56 & 53.36 & 56.32 & 64.02 & 59.8 & 59.36 & 49.13 & 69.3 & 53.41 \\  \multicolumn{2}{c}{DRL 4} & 49.08 & 54.34 & 57.16 & 65.26 & 61.41 & 59.35 & 55.23 & 46.19 & 39.59 \\  \multicolumn{2}{c}{DRL 5} & 34.31 & 41.62 & 49.56 & 56.95 & 60.46 & 58.74 & 49.18 & 68.99 & 53.01 \\  \multicolumn{2}{c}{DRL 6} & 33.85 & 40.31 & 48.7 & 59.89 & 63.3 & 59.03 & 55.66 & 46.57 & 40.0 \\ 
\bottomrule 
\end{tabular}
\end{table}

The results of the experiments are presented in Tables \ref{tab:app_ATDO}, \ref{tab:app_AOCT}, and \ref{tab:app_PUO}, corresponding to the Average Travel Distance per Order (ATDO), Average Order Completion Time (AOCT), and Percentage of Unfulfilled Orders (PUO) metrics, respectively. These results represent the averages over 10 instances for arrival rates $\lambda = 0.01, 0.02, \dots, 0.09$.
As anticipated, the experimental findings suggest that higher initial pick sizes tend to be restrictive, resulting in longer travel distances, increased order completion times, and lower order fulfillment ratios. Consequently, we recommend setting the initial pick size to 1 for this dynamic order-picking scenario.
An important question arises: why impose an initial pick size restriction at all? Allowing pickers unrestricted movement, even when no orders are queued, might appear beneficial. However, this could lead to unnecessary empty travels and cause training delays for the DRL agent, particularly under lower order arrival rates.
Regarding the performance of the models with initial pick size of unit value, DRL 5 (without cross-aisle rerouting) and DRL 6 (with cross-aisle rerouting) demonstrate comparable results on average. Nonetheless, at higher order arrival rates, the unrestricted variant (DRL 6) generally yields better performance.

\begin{table}[]
    \centering
    \scriptsize
    \caption{Sensitivity Analysis: Average Order Completion Time (AOCT)}
    \label{tab:app_AOCT}
    \begin{tabular}{cc rrrrrrrrr}
    \toprule
\multicolumn{2}{c}{}                & \multicolumn{9}{c}{AOCT (seconds) w.r.t. $\lambda$} \\
 \cmidrule(lr){3-11}
\multicolumn{2}{c}{Model} & \multicolumn{1}{c}{0.01}   & \multicolumn{1}{c}{0.02}   & \multicolumn{1}{c}{0.03}   & \multicolumn{1}{c}{0.04}   & \multicolumn{1}{c}{0.05}   & \multicolumn{1}{c}{0.06}   & \multicolumn{1}{c}{0.07}  & \multicolumn{1}{c}{0.08}  & \multicolumn{1}{c}{0.09}  \\
 \midrule
\multicolumn{2}{c}{DRL 1} & 480.7 & 300.8 & 254.2 & 254.9 & 226.3 & 282.2 & 295.8 & 505.0 & 581.6 \\  \multicolumn{2}{c}{DRL 2} & 474.5 & 304.5 & 252.6 & 254.9 & 233.4 & 257.2 & 299.1 & 366.6 & 538.2 \\  \multicolumn{2}{c}{DRL 3} & 333.8 & 230.1 & 201.7 & 213.9 & 216.8 & 292.3 & 307.5 & 510.0 & 587.7 \\  \multicolumn{2}{c}{DRL 4} & 329.5 & 227.6 & 200.9 & 214.7 & 221.8 & 263.8 & 321.9 & 367.5 & 557.7\\  \multicolumn{2}{c}{DRL 5} & 66.8 & 90.1 & 121.5 & 164.7 & 206.5 & 290.3 & 305.7 & 524.2 & 579.6 \\  \multicolumn{2}{c}{DRL 6} & 67.4 & 87.2 & 118.1 & 171.6 & 210.2 & 260.9 & 323.1 & 369.0 & 556.3 \\ 
\bottomrule 
\end{tabular}
\end{table}

\begin{table}[]
    \centering
    \scriptsize
    \caption{Sensitivity Analysis: Percentage of Unfulfilled Orders (PUO)}
    \label{tab:app_PUO}
    \begin{tabular}{cc rrrrrrrrr}
 \toprule
\multicolumn{2}{c}{}                & \multicolumn{9}{c}{PUO (\%) w.r.t. $\lambda$} \\
 \cmidrule(lr){3-11}
\multicolumn{2}{c}{Model} & \multicolumn{1}{c}{0.01}   & \multicolumn{1}{c}{0.02}   & \multicolumn{1}{c}{0.03}   & \multicolumn{1}{c}{0.04}   & \multicolumn{1}{c}{0.05}   & \multicolumn{1}{c}{0.06}   & \multicolumn{1}{c}{0.07}  & \multicolumn{1}{c}{0.08}  & \multicolumn{1}{c}{0.09}  \\
 \midrule
\multicolumn{2}{c}{DRL 1} & 1.86 & 1.08 & 0.95 & 1.07 & 0.78 & 1.02 & 0.97 & 1.73 & 2.16 \\  \multicolumn{2}{c}{DRL 2} & 1.80 & 1.14 & 0.84 & 0.93 & 0.76 & 0.68 & 0.90 & 1.31 & 2.00 \\  \multicolumn{2}{c}{DRL 3} & 1.37 & 0.93 & 0.98 & 0.66 & 0.81 & 0.94 & 0.98 & 1.76 & 2.08 \\  \multicolumn{2}{c}{DRL 4} & 1.33 & 0.91 & 0.92 & 0.73 & 0.74 & 0.74 & 0.84 & 1.36 & 2.19 \\  \multicolumn{2}{c}{DRL 5} & 0.36 & 0.31 & 0.43 & 0.75 & 0.69 & 0.86 & 1.03 & 1.71 & 2.01 \\  \multicolumn{2}{c}{DRL 6} & 0.64 & 0.36 & 0.30 & 0.60 & 0.56 & 0.69 & 1.01 & 1.25 & 2.27 \\ 
\bottomrule 
\end{tabular}
\end{table}

\section{Analysis of the DRL Model with Fluctuating Arrival Rates} 
\label{app:analysis_varying_rates}

In practical warehouse operations, demand rates often fluctuate throughout the operational horizon, characterized by peak and trough periods. Although our proposed DRL framework has been trained and rigorously benchmarked under specific order arrival rates, we also evaluate its performance in scenarios that more accurately reflect real-world dynamics. In this section, we specifically investigate the influence of varying order arrival rates on the proposed DRL model by simulating realistic operational conditions.

To effectively replicate practical variations in demand, we consider two cases with distinct arrival rate sequences: Case 1 ($\{0.02, 0.06, 0.04, 0.08\}$) and Case 2 ($\{0.09, 0.01, 0.03, 0.07\}$). Each arrival rate within a sequence is maintained for two hours, cumulatively representing a full 8-hour work shift. For example, in Case 1, the arrival rate is 0.02 for the first two hours, 0.06 for the next two, followed by 0.04 and 0.08 in subsequent two-hour intervals. During each 2-hour period, we utilize the DRL model previously trained on the arrival rate closest to the observed rate. For instance, when the arrival rate is either 0.01 or 0.02, the DRL model trained at an arrival rate of 0.02 is employed.
Tables \ref{tab:arrival_rate_seq1} and \ref{tab:arrival_rate_seq2} summarize the results of these experiments.

\begin{table}[htbp]
    \centering
    \caption{Performance metrics for different arrival rate sequences}
    \scriptsize
    \label{tab:performance_metrics}
    \begin{subtable}[t]{0.48\textwidth}
        \centering
        \caption{Arrival rate sequence \{0.02, 0.06, 0.04, 0.08\}}
        \label{tab:arrival_rate_seq1}
        \begin{tabular}{lccc}
            \toprule
            \multirow{2}{*}{Model} & ATDO & AOCT & PUO \\ 
                                   & (meters) & (seconds) & (\%) \\
            \midrule
            DRL        & 60.52 & 534.80  & 1.41 \\
            Baseline 1 & 8.19  & 1,581.31 & 10.26 \\
            Baseline 2 & 8.30  & 1,484.38 & 9.99 \\
            Baseline 3 & 8.29  & 1,485.79 & 9.99 \\
            Baseline 4 & 8.29  & 1,478.38 & 9.78 \\
            Baseline 5 & 8.29  & 1,479.00 & 9.78 \\
            \bottomrule
        \end{tabular}
    \end{subtable}%
    \hfill
    \begin{subtable}[t]{0.48\textwidth}
        \centering
        \caption{Arrival rate sequence \{0.09, 0.01, 0.03, 0.07\}}
        \scriptsize
        \label{tab:arrival_rate_seq2}
        \begin{tabular}{lccc}
            \toprule
            \multirow{2}{*}{Model} & ATDO & AOCT & PUO \\ 
                                   & (meters) & (seconds) & (\%) \\
            \midrule
            DRL        & 68.50 & 415.10  & 0.91 \\
            Baseline 1 & 8.21  & 565.71  & 1.42 \\
            Baseline 2 & 8.82  & 467.51  & 1.22 \\
            Baseline 3 & 8.82  & 462.46  & 1.22 \\
            Baseline 4 & 8.86  & 466.82  & 1.24 \\
            Baseline 5 & 8.86  & 464.02  & 1.21 \\
            \bottomrule
        \end{tabular}
    \end{subtable}
\end{table}

The experimental results demonstrate the adaptability of the DRL models under dynamically changing order arrival conditions, highlighting their robustness compared to traditional baseline models. In Case 1, where the arrival rate sequence follows \{0.02, 0.06, 0.04, 0.08\}, the baseline models perform poorly, with significantly higher throughput times and a much larger percentage of unfulfilled orders compared to the DRL model. The DRL model achieves a 66.2\% reduction in throughput time and a 85.6\% decrease in unfulfilled orders relative to the best-performing baseline.
In Case 2, where the arrival rate sequence is \{0.09, 0.01, 0.03, 0.07\}, the baseline models exhibit relatively improved performance. This can be attributed to the highest arrival rate occurring in the initial phase, leading to a backlog of unfulfilled orders that can be processed during the subsequent phases when order arrivals decrease. However, even in this scenario, the DRL model continues to outperform the baselines, achieving a 10.2\% reduction in throughput time and a 24.8\% decrease in unfulfilled orders compared to the best-performing baseline.

These findings underscore the effectiveness of the proposed DRL approach in handling fluctuating demand patterns, ensuring lower order fulfillment delays and improved operational efficiency across different demand scenarios.
\color{black}
\end{document}